\documentclass[a4paper,11pt,twoside]{amsart}
\usepackage[english]{babel}
\usepackage[utf8]{inputenc}

\usepackage[a4paper,inner=3cm,outer=3cm,top=4cm,bottom=4cm,pdftex]{geometry}
\usepackage{fancyhdr}
\usepackage{color}
\usepackage{bold-extra}
\usepackage{mathrsfs}
\usepackage{comment}
\usepackage{graphics}
\usepackage{aliascnt}

\usepackage{amsmath, amsfonts, amssymb, amsthm}
\usepackage{mathtools}
\usepackage{enumerate}
\usepackage[shortlabels]{enumitem}
\usepackage{cite}
\usepackage[pdftex,citecolor=green,linkcolor=red]{hyperref}
\def\Z{\mathbb Z}

\def\N{\mathbb N}

\newcommand{\Primes}{\mathscr P}
\def\1{{\bf 1}}

\theoremstyle{plain}
\newtheorem{theorem}{Theorem}
\newtheorem{proposition}{Proposition}
\newtheorem{lemma}{Lemma}
\newtheorem{corollary}[theorem]{Corollary}
\theoremstyle{definition}
\newtheorem*{acknowledgment}{Acknowledgments}
\theoremstyle{remark}

\begin{document}
\title[The critical-window profile for $d_k$ in short intervals]
{The critical-window profile for the $k$-fold divisor function in almost all short intervals}

\author{Yu-Chen Sun}
\address {Department of Mathematics and Statistics, University of Turku, 20014 Turku, Finland}
\email{yuchensun93@163.com}


\begin{abstract}
We establish an almost-all transition theorem for the $k$-fold divisor function
in short intervals.  Let $k\geq2$ be fixed, let
$\ell=\log\log X$, and set
\[
        D_k(X)=(\log X)^{k\log k-k+1},
        \qquad
        M_k(x)=\frac1x\sum_{x<n\leq2x}d_k(n).
\]
The critical scale is $D_k(X)$.  In the bounded part of the transition window,
put
\[
        A(X)=\frac1{\sqrt{\ell}}\log\frac{h}{D_k(X)}
\]
and assume that $A(X)=O(1)$.  Then, for almost all integers $x\in[X,2X]$,
\[
        \frac1h\sum_{x<n\leq x+h}d_k(n)
        =
        \left(\Phi_{\rm G}\left(\frac{A(X)}{\sqrt{k}\log k}\right)+o(1)\right)
        M_k(x),
\]
where $\Phi_{\rm G}$ denotes the standard Gaussian distribution function.  Above
the window, namely when
\[
        \frac{\log(h/D_k(X))}{\sqrt{\ell}}\to+\infty,
\]
the full long average is recovered for almost all $x$.
\end{abstract}
\maketitle

\section{Introduction}\label{intro}
\setcounter{lemma}{0}\setcounter{theorem}{0}\setcounter{proposition}{0}\setcounter{equation}{0}

A central problem in analytic number theory is to compare local averages of
arithmetic functions with their long averages.  The celebrated theorem of Matom\"aki and
Radziwi{\l\l} \cite[Theorem 1]{MR16} shows that bounded multiplicative functions
have the expected average on almost all short intervals whose lengths tend to
infinity.  For unbounded functions such as the divisor functions, one expects a
new threshold phenomenon, already at polylogarithmic scales.

The divisor functions are among the basic unbounded multiplicative functions.
For fixed $k\geq2$ let
\[
        d_k(n)=\#\{(m_1,\ldots,m_k)\in\mathbb N^k:m_1\cdots m_k=n\};
\]
equivalently, $d_k$ is the coefficient sequence of $\zeta(s)^k$.  Short
interval averages of $d_k$ measure how far the local distribution of
factorisations can be predicted from the corresponding long average.  The
difficulty is that $d_k$ has a heavy tail: integers with many prime factors
carry disproportionately large weight.  This makes the bounded-function theory
insufficient by itself and leads to the critical scale studied below.

In this paper we determine the almost-all transition scale for the $k$-fold
divisor function.  All logarithms are natural.  Throughout the paper
$X\to\infty$ and the integer $k\geq2$ is fixed.  We write $n\sim N$ for
$N<n\leq2N$.  Let $\Primes$ denote the set of primes, and let
\[
        \omega(n)=\#\{p\in\Primes:p\mid n\},
        \qquad
        \Omega(n)=\sum_{p^a\Vert n}a .
\]
Thus $\omega(n)$ and $\Omega(n)$ count prime factors without and with
multiplicity, respectively.

We study how short an interval can be while still recovering the long average of
$d_k$.  More precisely, for almost all $x\in[X,2X]$ we ask when
\begin{equation}\label{eq:short_interval_question}
        \frac1h\sum_{x<n\leq x+h}d_k(n)
        \sim
        \frac{1}{x}\sum_{n\sim x}d_k(n)
        =: M(x),
\end{equation}
where $M(x)\sim P_{k-1}(\log x)$, with $P_{k-1}$ a polynomial of degree $k-1$.

Polynomial-scale analogues for divisor-type functions are known.  Closer to the
present paper, Mangerel \cite{mangerel2021divisor} extended the
Matom\"aki--Radziwi{\l\l} theorem to divisor-bounded multiplicative functions.
He proved that \eqref{eq:short_interval_question} holds
for almost all $x\in[X,2X]$ whenever
\[
        \frac{h}{(\log X)^{(k-1)^2}}\to+\infty .
\]
Our endpoint theorem replaces this length by the condition
\[
        \frac{\log(h/D_k(X))}{\sqrt{\log\log X}}\to+\infty .
\]
At the lower edge this is
$h=(\log X)^{k\log k-k+1+o(1)}$.  Since
\[
        (k-1)^2-(k\log k-k+1)=k(k-1-\log k)>0
\]
for every $k\geq2$, this gives a power-saving improvement over Mangerel's
exponent for the divisor function.  Theorem \ref{critical_window} goes further
by identifying the transition profile.  Matom\"aki, Shao, Tao and Ter\"av\"ainen
\cite[Theorem 1.1(iii)]{MSTT23} prove all-interval higher-uniformity results
for $d_k-d_k^\sharp$, and the almost-all version is proved in
\cite[Theorem 1.1(v)]{MRSTT26}.  These results do not address the
polylogarithmic scale at which the unboundedness of $d_k$ creates a transition.

One technical point is that the natural $\Omega$-cutoff cannot be treated as a
small tail in the critical range.  We keep the cutoff through the Ramar\'e
decomposition and open it only on the final residual sum; after this reduction,
the remaining large-prime layer gives a uniform high-frequency Hal\'asz saving.

The goal of the present paper is to identify this transition exactly.  In the
statements below
\[
        \ell=\log\log X,\qquad
        D_k(X)=(\log X)^{k\log k-k+1}.
\]
The length $h=h(X)$ is always a positive integer.  For real $X\leq x\leq2X$ let
\[
        M_k(x)=\frac1x\sum_{x<n\leq2x}d_k(n).
\]
Whenever an
exceptional set $E_X$ is specified, the corresponding error term is uniform for
$x\in([X,2X]\cap\mathbb Z)\setminus E_X$.  The set $E_X$ is allowed to change
from one occurrence to the next, but in each occurrence it satisfies
$|E_X|=o(X)$.

Write
\[
        \Phi_{\rm G}(t)=\frac1{\sqrt{2\pi}}\int_{-\infty}^{t}\exp(-u^2/2)\,du
\]
for the standard Gaussian distribution function.

\begin{theorem}[Critical-window profile]\label{critical_window}
Let $h=h(X)\in\mathbb N$ satisfy $1\leq h\leq X$, and let
\[
        A(X)=\frac1{\sqrt{\ell}}\log\frac{h}{D_k(X)}.
\]
Assume that $A(X)=O(1)$.
Then there exists a set $E_X\subset [X,2X]\cap\mathbb Z$ with $|E_X|=o(X)$ such
that, uniformly for $x\in([X,2X]\cap\mathbb Z)\setminus E_X$,
\[
        \frac1h\sum_{x<n\leq x+h}d_k(n)
        =
        \left(\Phi_{\rm G}\left(\frac{A(X)}{\sqrt{k}\log k}\right)+o(1)\right)
        M_k(x).
\]
\end{theorem}

The endpoint regime corresponds formally to letting $A(X)$ tend to
$+\infty$, where the Gaussian factor tends to $1$.  The following corollary
proves this limiting statement uniformly.

\begin{corollary}\label{dvthm}
Let $h=h(X)\in\mathbb N$ satisfy
\[
        1\leq h\leq X,
        \qquad
        \frac{\log(h/D_k(X))}{\sqrt{\ell}}\to+\infty .
\]
Then
\[
\frac1X\int_X^{2X}\left|
\frac1h\sum_{x<n\leq x+h}d_k(n)-M_k(x)
\right|dx=o((\log X)^{k-1}).
\]
Moreover, there exists a set $E_X\subset [X,2X]\cap\mathbb Z$ with
$|E_X|=o(X)$ such that, uniformly for
$x\in([X,2X]\cap\mathbb Z)\setminus E_X$,
\[
\frac1h\sum_{x<n\leq x+h}d_k(n)-M_k(x)
=o((\log X)^{k-1}).
\]
\end{corollary}

Theorem \ref{critical_window} describes the finite part of the transition around
the critical scale $D_k(X)$.  The endpoint estimate in Corollary \ref{dvthm}
corresponds to the limiting regime
\[
        \frac{\log(h/D_k(X))}{\sqrt{\ell}}\to+\infty,
\]
where the Gaussian factor tends to $1$.

There are two key ingredients in the proof: a weighted Erd\H{o}s--Kac theorem
for $\Omega(n)$ under the $d_k$-measure, and an $\Omega$-truncated
Matom\"aki--Radziwi{\l\l} estimate.  These are explained in Section
\ref{sec:outline}; the main results are then deduced in Section
\ref{structure}.

\begin{acknowledgment}
The author is grateful to his supervisor Kaisa Matom\"aki for suggesting this project, drawing his attention to Matom\"aki--Radziwi{\l\l}--Tao \cite{MR3961326}, and many useful discussions. He also thanks Joni Ter\"av\"ainen and Martin \v Cech for their useful comments. During the work the author was supported by UTUGS funding, working in the Academy of Finland project no. 333707.
\end{acknowledgment}

\section{Proof ingredients and outline}\label{sec:outline}
\setcounter{lemma}{0}\setcounter{theorem}{0}\setcounter{proposition}{0}\setcounter{equation}{0}

We first explain the source of the Gaussian factor.  Consider the probability
measure on integers $n\in[X,2X]$ defined by
\[
        \mathbb P_{k,X}(n)
        =
        \frac{d_k(n)}{\sum_{X<m\leq2X}d_k(m)}.
\]
Under this $d_k$-weighted measure, the random variable $\Omega(n)$ satisfies the
central limit theorem used below:
\[
        \frac{\Omega(n)-k\log\log X}{\sqrt{k\log\log X}}
        \Longrightarrow N(0,1).
\]
This is the distributional input recorded in Lemma \ref{weighted_EK}(i) below,
which is a weighted Erd\H{o}s--Kac theorem for $d_k$.

The second ingredient is analytic.  If we truncate to
\[
        B(L)=\{n\in\mathbb N:\Omega(n)\leq L\},
\]
then $d_k(n)\leq k^L$ on $B(L)$.  For the corresponding Dirichlet polynomial
\[
        F(s)=\sum_{n\sim X}\frac{d_k(n)\1_{\Omega(n)\le L}}{n^s},
\]
Lemma \ref{mean_value_B}, applied with $M,Y\asymp X$, gives the model bound
\[
        \int_{-T}^{T}|F(1+it)|^2dt
        \ll_k
        \left(\frac{T}{X}H+1\right)(\log X)^{2k-2},
\]
where
\[
        H=H(L;X):=
        \frac{1}{X(\log X)^{2k-2}}
        \sum_{n\sim X}d_k(n)^2\1_{\Omega(n)\le L}
        \ll_k \frac{k^L}{(\log X)^{k-1}}.
\]
For a short interval of length $h$, Perron's formula leads to heights
$T\asymp X/h$.  Thus the first term in the mean-value bound is acceptable when
\[
        \frac{TH}{X}\asymp \frac{H}{h}\ll1.
\]
It is therefore enough that $h$ exceeds the upper-bound length.  In the case
$L=k\ell+u\sqrt{\ell}$ this pushes $h$ beyond
\[
        \frac{k^L}{(\log X)^{k-1}}
        =D_k(X)\exp(u\log k\sqrt{\ell}).
\]
Consequently, if the interval length is
\[
        h=D_k(X)\exp(A\sqrt{\ell}+o(\sqrt{\ell}))
\]
then the condition $h\gg k^L/(\log X)^{k-1}$ forces the admissible parameter to
satisfy $u<A/\log k$, up to the small margin used in the proof.  Thus one can
control, in the almost-all sense, the part of the $d_k$-mass for which
\[
        \Omega(n)\leq k\ell+\frac{A}{\log k}\sqrt{\ell},
\]
up to a negligible transition band.  The weighted Erd\H{o}s--Kac theorem then
gives
\[
        \mathbb P_{k,X}\left(
        \Omega(n)\leq k\ell+\frac{A}{\log k}\sqrt{\ell}
        \right)
        =
        \Phi_{\rm G}\left(\frac{A}{\sqrt{k}\log k}\right)+o(1),
\]
which is the factor in Theorem \ref{critical_window}.

The rigorous argument follows this outline.  For the distributional input we use
Khan, Milinovich and Subedi
\cite[Corollary 1.2, formula (1.5)]{KMS22}; in Lemma \ref{weighted_EK} we check
that their Selberg--Delange input is unchanged when $\omega(n)$ is replaced by
$\Omega(n)$.  The
analytic input is proved by Perron's formula and a Ramar\'e decomposition.  The
latter remains compatible with the cutoff because, for $p\nmid m$,
\[
        d_k(mp)=k\,d_k(m),\qquad \Omega(mp)=\Omega(m)+1.
\]
The main new technical point, isolated in Lemma
\ref{final_residual_pointwise}, concerns the final residual.  A direct
estimate for
\[
        \frac1X\sum_{\substack{n\sim X\\ \Omega(n)>L}}d_k(n)
\]
does not give a logarithmic saving in the central range of $L$.  Thus the
$\Omega$-cutoff cannot simply be discarded as a negligible tail.  Instead, only
after the final Ramar\'e decomposition, we apply Cauchy's formula to the
residual condition $\Omega(m)\le L-1$ and write it as an integral of the
multiplicative weights $z^{-\Omega(m)}$.  Each resulting integrand is a
divisor-bounded multiplicative function whose prime layer still contains
primes up to length $X^{1-o(1)}$.  This large-prime layer forces a uniform
high-frequency Hal\'asz distance, even after the local prime-cell restrictions
are imposed.  This is the mechanism that supplies the logarithmic saving in the
final large-value range.
Section \ref{structure} deduces the main results from these inputs, and Section
\ref{pf_MR_bd} proves the required mean-square estimate.

\section{Deduction of the main results from the mean-square estimates}\label{structure}
\setcounter{lemma}{0}\setcounter{theorem}{0}\setcounter{proposition}{0}\setcounter{equation}{0}

In this section we prove Theorem \ref{critical_window}.  The analytic inputs are the
truncated mean-square estimate stated below as Proposition
\ref{MR_bound_flexible}; it is proved in Section \ref{pf_MR_bd}.

Throughout the rest of the paper we retain the notation and conventions from the
introduction.  In particular $k\geq2$ is fixed,
\[
        \ell=\log\log X,\qquad
        D_k(X)=(\log X)^{k\log k-k+1},
        \qquad
        M_k(x)=\frac1x\sum_{x<n\leq2x}d_k(n).
\]
All short lengths $h$ are positive integers depending on $X$.  We fix
\begin{equation}\label{sigma_def}
        \sigma=\frac1{100},
        \qquad
        \beta_0=\frac1{20},
        \qquad
        h_1=\lfloor X^{3/4}\rfloor,
        \qquad
        \delta_0=X^{-\beta_0}.
\end{equation}
We shall use the following two consequences of Shiu's theorem.

\begin{lemma}[Shiu-type estimates]\label{Shiu_bound}
Let $0\leq \eta\leq1$, $\sqrt Y<y\leq Y$, and $Y\asymp X$. Then
\begin{equation}\label{Shiu_weighted}
        \sum_{Y-y<n\leq Y}d_k(n)(1+\eta)^{\Omega(n)}
        \ll_k y(\log X)^{k(1+\eta)-1}.
\end{equation}
Moreover, if $2\leq P\leq Q\leq X$, then
\begin{equation}\label{Shiu_no_prime}
        \sum_{\substack{X<n\leq 3X\\ p\mid n\Rightarrow p\notin [P,Q]}}d_k(n)
        \ll_k X(\log X)^{k-1}\left(\frac{\log P}{\log Q}\right)^k.
\end{equation}
\end{lemma}

\begin{proof}
Both estimates are direct applications of Shiu's theorem
\cite[Theorem 1]{Sh80}.  For \eqref{Shiu_weighted} use
$f(n)=d_k(n)(1+\eta)^{\Omega(n)}$, whose prime factor is
$1+k(1+\eta)/p+O_k(p^{-2})$.  For \eqref{Shiu_no_prime} use
$f(n)=d_k(n)\1_{(n,\prod_{P\leq p\leq Q}p)=1}$; Mertens' formula gives
\[
        \prod_{P\leq p\leq Q}\left(1-\frac{k}{p}+O_k\left(\frac1{p^2}\right)\right)
        \ll_k\left(\frac{\log P}{\log Q}\right)^k
\]
\cite[Theorem I.1.6, p.~17]{Te15}.
\end{proof}

We next record the distributional estimates for $\omega(n)$ and $\Omega(n)$
under the $d_k$-weight.

\begin{lemma}[Weighted distribution estimates]
\label{weighted_EK}\label{large_omega_count}
The following estimates hold.
\begin{enumerate}[(i)]
\item If $z(n)\in\{\omega(n),\Omega(n)\}$, then for every fixed $a<b$,
\[
\frac{\sum_{X<n\leq2X}d_k(n)\1_{a\sqrt{k\ell}<z(n)-k\ell\leq b\sqrt{k\ell}}}
{\sum_{X<n\leq2X}d_k(n)}
\to \Phi_{\rm G}(b)-\Phi_{\rm G}(a).
\]
In particular, if $M=M(X)\to\infty$, then
\[
\frac{\sum_{X<n\leq2X}d_k(n)\1_{|z(n)-k\ell|\leq M\sqrt{\ell}}}
{\sum_{X<n\leq2X}d_k(n)}
\to1.
\]
\item If $z(n)\in\{\omega(n),\Omega(n)\}$, then uniformly for $|u|\leq \ell^{1/3}$,
\[
        \#\{X<n\leq3X:z(n)\geq k\ell+u\sqrt{\ell}\}
        \ll_k
        \frac{X}{D_k(X)\exp(u\log k\sqrt{\ell})}.
\]
\end{enumerate}
\end{lemma}

\begin{proof}
Part (i) for $z(n)=\omega(n)$ is precisely the weighted Erd\H{o}s--Kac theorem
of Khan, Milinovich and Subedi
\cite[Corollary 1.2, formula (1.5)]{KMS22}, applied at $2X$ and $X$ and
subtracted.  For $z(n)=\Omega(n)$, it is enough to record that the
Selberg--Delange input used in their proof has the same form.  For $w$ in a
fixed neighbourhood of $1$,
\[
        \sum_{n=1}^{\infty}\frac{d_k(n)w^{\Omega(n)}}{n^s}
        =
        \prod_p(1-wp^{-s})^{-k}.
\]
We factor this product as
\[
        \prod_p(1-wp^{-s})^{-k}
        =
        \zeta(s)^{kw}\zeta(2s)^{k(w^2-w)/2}\widetilde G_k(s,w).
\]
After taking logarithms, the factors $\zeta(s)^{kw}$ and
$\zeta(2s)^{k(w^2-w)/2}$ remove the $p^{-s}$ and $p^{-2s}$ terms, respectively.
Thus
\[
        \log \widetilde G_k(s,w)=\sum_p O_{k,w}(p^{-3s}),
\]
uniformly for $w$ in compact subsets of this neighbourhood.  Hence
$\widetilde G_k(s,w)$ is holomorphic and locally bounded in a half-plane
$\Re s>1/3+\varepsilon$.  This is exactly the Selberg--Delange structure
required in the proof of \cite[Corollary 1.2, formula (1.5)]{KMS22}, and their
characteristic-function argument applies without change to
$d_k(n)w^{\Omega(n)}$.

For part (ii) we use the Sathe--Selberg upper bound
\[
        \#\{X<n\le 2X:z(n)=m\}
        \ll_C
        \frac{X}{\log X}\frac{\ell^{m-1}}{(m-1)!},
        \qquad 1\le m\le C\ell,
\]
for $z(n)\in\{\omega(n),\Omega(n)\}$, which follows from
\cite[Chapter~II.6, Theorems 4 and 5]{Te15}.  Let
$m_0=\lceil k\ell+u\sqrt\ell\rceil$.  Since $|u|\le \ell^{1/3}$ and $k\ge2$,
we have, uniformly for $m\ge m_0$ with $m\le C\ell$,
\[
        \frac{\ell^m/m!}{\ell^{m-1}/(m-1)!}
        =\frac{\ell}{m}\le \frac1{k+o(1)}<1.
\]
Thus the tail up to $C\ell$ is bounded by the first term.  Stirling's formula
then gives
\[
        \sum_{m_0\le m\le C\ell}
        \frac{X}{\log X}\frac{\ell^{m-1}}{(m-1)!}
        \ll_k
        X\exp\{-(k\log k-k+1)\ell-u\log k\sqrt{\ell}\}.
\]
The range $m>C\ell$ is harmless if $C=C(k)$ is large.  Indeed, for fixed
$R>1$, Rankin's trick and Shiu's theorem give
\[
        \#\{X<n\le2X:z(n)>C\ell\}
        \le R^{-C\ell}\sum_{X<n\le2X}R^{z(n)}
        \ll_R X(\log X)^{R-1-C\log R},
\]
which is smaller than the displayed bound after increasing $C$.  Applying the
same argument with $X$ replaced by $2X$ covers $[2X,3X]$, since
\(\log\log(2X)=\ell+O(1/\log X)\).
\end{proof}

\begin{lemma}[Thin weighted bands]\label{thin_weighted_band}
Let $B_0>0$, let $\eta(X)\to0$, and let $|\nu(X)|\le B_0$.  Then
\[
        \sum_{\substack{X<n\le3X\\
        |\Omega(n)-k\ell-\nu(X)\sqrt\ell|\le \eta(X)\sqrt\ell}}
        d_k(n)
        =
        o\bigl(X(\log X)^{k-1}\bigr).
\]
\end{lemma}

\begin{proof}
It suffices to treat $X<n\le2X$, the interval $2X<n\le3X$ following by the
same argument with $X$ replaced by $2X$.  Let
\[
        F_X(t)=
        \frac{\sum_{X<n\le2X}d_k(n)
        \1_{(\Omega(n)-k\ell)/\sqrt{k\ell}\le t}}
        {\sum_{X<n\le2X}d_k(n)} .
\]
Part (i) of Lemma \ref{weighted_EK}, together with Polya's theorem, gives
$F_X(t)\to\Phi_{\rm G}(t)$ uniformly for $t$ in compact intervals.  Put
\[
        a_X=\frac{\nu(X)-\eta(X)}{\sqrt k},\qquad
        b_X=\frac{\nu(X)+\eta(X)}{\sqrt k}.
\]
For any fixed $\delta>0$ and all large $X$,
\[
        F_X(b_X)-F_X(a_X-\delta)
        \le
        \sup_{|t|\le B_0+1}|F_X(t)-\Phi_{\rm G}(t)|
        +\Phi_{\rm G}(b_X)-\Phi_{\rm G}(a_X-\delta).
\]
Taking first $X\to\infty$ and then $\delta\downarrow0$ gives the claim on
$[X,2X]$.  For $[2X,3X]$ one uses
$\ell_2=\log\log(2X)=\ell+o(1)$ and rewrites the centre as
$k\ell_2+\nu_2\sqrt{\ell_2}$ with $|\nu_2|\le B_0+1$ and band width
$(\eta+o(1))\sqrt{\ell_2}$.
\end{proof}

\begin{lemma}[Long averages with an $\Omega$-cutoff]\label{long_average_Omega_cutoff}
Let $\theta>7/12$ be fixed.  Let
\[
        L=k\ell+\nu(X)\sqrt{\ell},
        \qquad |\nu(X)|\leq \ell^{1/3},
\]
and let $B=\{n\in\mathbb N:\Omega(n)\leq L\}$.  Uniformly for
$X\leq x\leq2X$ and
\[
        X^\theta\leq y\leq X,
\]
one has
\[
        \frac1y\sum_{\substack{x<n\leq x+y\\ n\in B}}d_k(n)
        =\left(\Phi_{\rm G}\left(\frac{\nu(X)}{\sqrt k}\right)+o(1)\right)
        M_k(x).
\]
\end{lemma}

\begin{proof}
We use the Selberg--Delange input quoted as Lemma 2.1 of Liu and Wu
\cite{LiuWu21}, together with their smoothing argument in Section~4.  Applied
to
\[
        \sum_{n=1}^{\infty}\frac{d_k(n)z^{\Omega(n)}}{n^s}
        =
        \prod_p(1-zp^{-s})^{-k},
\]
that input gives, uniformly for $X\le x\le2X$, $X^\theta\le y\le X$, and $z$
in a fixed neighbourhood of $1$,
\begin{equation}\label{short_SD_Omega}
        \sum_{x<n\le x+y}d_k(n)z^{\Omega(n)}
        =
        y(\log x)^{kz-1}
        \left(C_k(z)+O_k((\log X)^{-c_\theta})\right),
\end{equation}
where $c_\theta>0$ and
\[
        C_k(z)=
        \frac{\zeta(2)^{k(z^2-z)/2}\widetilde G_k(1,z)}{\Gamma(kz)}.
\]
Here $\widetilde G_k$ is defined by
\[
        \prod_p(1-zp^{-s})^{-k}
        =
        \zeta(s)^{kz}\zeta(2s)^{k(z^2-z)/2}\widetilde G_k(s,z).
\]
The logarithm of $\widetilde G_k(s,z)$ starts with terms $O_z(p^{-3s})$; hence
$\widetilde G_k(s,z)$ is holomorphic and uniformly bounded in the region
required for Liu--Wu's Lemma 2.1, and $C_k(1)=1/\Gamma(k)$.

Let
\[
        F_{X,x,y}(t)=
        \frac{\sum_{x<n\le x+y}d_k(n)
        \1_{(\Omega(n)-k\ell)/\sqrt{k\ell}\le t}}
        {\sum_{x<n\le x+y}d_k(n)} .
\]
For bounded $\tau$ put $z=\exp(i\tau/\sqrt{k\ell})$.  Dividing
\eqref{short_SD_Omega} by its value at $z=1$ gives the characteristic function
of $(\Omega(n)-k\ell)/\sqrt{k\ell}$ under the $d_k$-weight on $(x,x+y]$:
\[
        e^{-i\tau\sqrt{k\ell}}\,
        \frac{\sum_{x<n\le x+y}d_k(n)z^{\Omega(n)}}
             {\sum_{x<n\le x+y}d_k(n)}
        =
        \exp(-\tau^2/2)+o(1),
\]
uniformly for $X\le x\le2X$ and $X^\theta\le y\le X$.  Here we used
$\log\log x=\ell+o(1)$ and $C_k(z)/C_k(1)=1+o(1)$.  The Berry--Esseen
smoothing argument of \cite[Section~4]{LiuWu21} therefore yields
\begin{equation}\label{short_distribution_uniform}
        \sup_{\substack{X\le x\le2X\\ X^\theta\le y\le X}}
        \sup_{t\in\mathbb R}
        |F_{X,x,y}(t)-\Phi_{\rm G}(t)|\to0.
\end{equation}
Taking $t=\nu(X)/\sqrt k$ in \eqref{short_distribution_uniform} gives the
cutoff asymptotic, uniformly for $|\nu(X)|\le\ell^{1/3}$.  Finally,
\eqref{short_SD_Omega} with $z=1$ gives
\[
        \sum_{x<n\le x+y}d_k(n)
        =
        y\frac{(\log x)^{k-1}}{\Gamma(k)}+o(y(\log X)^{k-1})
        =
        yM_k(x)+o(y(\log X)^{k-1})
\]
uniformly in the same range.
\end{proof}

The following Matom\"aki--Radziwi{\l\l} type estimates are proved in Section
\ref{pf_MR_bd}.  Since $h$ and $h_1$ are integers, their $L^1$ forms imply the
corresponding integer almost-all statements by Markov's inequality.  Since
$\Omega(n)$ is integer-valued, the condition $\Omega(n)\leq L$ is understood
literally for real $L$.

Let
\[
        P_1=\exp(\ell^{3/10}),\qquad
        Q_1=\exp\left(\frac{\ell^{2/5}}{1000}\right),
\]
and
\[
        P_2=\exp(\exp(\ell^{3/20})),\qquad
        Q_2=\exp(\exp(2\ell^{3/20})).
\]
We also set
\[
        P_3=\exp\left(\frac{\log X}{\ell^2}\right),\qquad
        Q_3=\exp\left(\frac{\log X}{100k^2\ell}\right),
\]
and denote by $\mathcal A$ the set of integers having at least one prime factor
in each of the intervals $[P_1,Q_1]$, $[P_2,Q_2]$, and $[P_3,Q_3]$.

\begin{proposition}[$\Omega$-truncated Matom\"aki--Radziwi{\l\l} estimate]
\label{MR_bound_flexible}
Let the prime blocks and $\mathcal A$ be as above.  Let
\[
        L=k\ell+\nu(X)\sqrt{\ell},\qquad |\nu(X)|\le \ell^{1/3},
\]
and define
\[
        B=\{n\in\N:\Omega(n)\leq L\},
        \qquad
        \mathcal S=B\cap\mathcal A,
        \qquad
        H=\frac{k^L}{(\log X)^{k-1}}.
\]
Let $h=h(X)\in\mathbb N$ satisfy
\[
        H\exp(\ell^{2/5})\leq h\leq h_1 .
\]
Then
\[
\frac1X\int_X^{2X}\left|
        \frac1h\sum_{\substack{x<n\leq x+h\\ n\in \mathcal S}}d_k(n)
        -
        \frac1{h_1}\sum_{\substack{x<n\leq x+h_1\\ n\in \mathcal S}}d_k(n)
\right|dx
=o((\log X)^{k-1}).
\]
Moreover, there exists a set $E_X\subset [X,2X]\cap\mathbb Z$ with
$|E_X|=o(X)$ such that, uniformly for
$x\in([X,2X]\cap\mathbb Z)\setminus E_X$,
\[
        \frac1h\sum_{\substack{x<n\leq x+h\\ n\in \mathcal S}}d_k(n)
        =
        \frac1{h_1}\sum_{\substack{x<n\leq x+h_1\\ n\in \mathcal S}}d_k(n)
        +o((\log X)^{k-1}).
\]
\end{proposition}

\begin{proof}[Proof of Theorem \ref{critical_window}, assuming Proposition \ref{MR_bound_flexible}]
Since $k\log k-k+1>0$ for $k\ge2$, the hypothesis $A(X)=O(1)$ gives
\[
        h=D_k(X)\exp(A(X)\sqrt\ell)
        =(\log X)^{k\log k-k+1}\exp(O(\sqrt\ell))\to\infty
\]
and $h=X^{o(1)}<h_1$ for all sufficiently large $X$.
Let
\[
        \eta=\frac{2\ell^{-1/10}}{\log k},
        \qquad
        \nu_\pm=\frac{A(X)}{\log k}\pm\eta,
        \qquad
        L_\pm=k\ell+\nu_\pm\sqrt{\ell},
\]
and let
\[
        B_\pm=\{n\in\N:\Omega(n)\leq L_\pm\},\qquad
        \mathcal S=B_-\cap\mathcal A,\qquad
        H_-=\frac{k^{L_-}}{(\log X)^{k-1}}.
\]
Then
\[
        H_-=D_k(X)\exp\{(A(X)-\eta\log k)\sqrt{\ell}\}
        =h\exp(-2\ell^{2/5}).
\]
Since $A(X)=O(1)$, for all sufficiently large $X$ one has
$|\nu_-|\le\ell^{1/3}$ and $H_-\exp(\ell^{2/5})\le h\le h_1$.  Hence
Proposition \ref{MR_bound_flexible} may be applied with $L=L_-$ and $B=B_-$.
By \eqref{Shiu_no_prime},
\begin{align*}
\sum_{\substack{X<n\leq3X\\ n\notin\mathcal A}}d_k(n)
&\ll X(\log X)^{k-1}
\left\{
\left(\frac{\log P_1}{\log Q_1}\right)^k+
\left(\frac{\log P_2}{\log Q_2}\right)^k+
\left(\frac{\log P_3}{\log Q_3}\right)^k
\right\}\\
&=o(X(\log X)^{k-1}).
\end{align*}
For $y\in\{h,h_1\}$ put
\[
        M_y(x)=
        \frac1y\sum_{\substack{x<n\le x+y\\ n\notin\mathcal A}}d_k(n).
\]
Then
\[
        \sum_{X\le x\le2X}M_y(x)
        \le
        \frac{y+1}{y}
        \sum_{\substack{X<n\le3X\\ n\notin\mathcal A}}d_k(n)
        =
        o\bigl(X(\log X)^{k-1}\bigr).
\]
Choosing any $\varepsilon_X\to0$ slowly enough and applying Markov's inequality
gives $M_y(x)=o((\log X)^{k-1})$ for both $y=h$ and $y=h_1$, outside
$o(X)$ integers.  Thus $B_-$ may be replaced by $\mathcal S$ in the two short
averages.  Proposition
\ref{MR_bound_flexible} therefore gives, for all but
$o(X)$ integers $x\in[X,2X]$,
\begin{equation}\label{critical_low_short}
        \frac1h\sum_{\substack{x<n\leq x+h\\ n\in B_-}}d_k(n)
        =
        \frac1{h_1}\sum_{\substack{x<n\leq x+h_1\\ n\in B_-}}d_k(n)+o((\log X)^{k-1}).
\end{equation}
Since $h_1\ge X^{2/3}$ for all large
$X$, Lemma \ref{long_average_Omega_cutoff} applies to the right hand side, for
example with $\theta=2/3$, and gives
\begin{equation}\label{critical_low_asymp}
        \frac1{h_1}\sum_{\substack{x<n\leq x+h_1\\ n\in B_-}}d_k(n)
        =
        \left(\Phi_{\rm G}\left(\frac{A(X)}{\sqrt k\log k}\right)+o(1)\right)
        M_k(x),
\end{equation}
because $\eta\to0$.  This gives the required lower bound for the full short average, since all terms are non-negative.
Also $M_k(x)\asymp_k(\log X)^{k-1}$ uniformly for $X\le x\le2X$, so every
additive $o((\log X)^{k-1})$ term below is $o(M_k(x))$.

For the upper bound let
\[
        R_\Omega=
        \{n\in[X,3X]:\Omega(n)>L_+\}.
\]
By Lemma \ref{large_omega_count}, with $z(n)=\Omega(n)$,
\[
        h\#R_\Omega
        \ll X\exp(-\eta\log k\sqrt{\ell}+O(1))=o(X),
\]
because $\eta\log k\sqrt\ell=2\ell^{2/5}$.  Let
\[
        \mathcal E_R=
        \{x\in[X,2X]\cap\Z:(x,x+h]\cap R_\Omega\ne\varnothing\}.
\]
Each $n\in R_\Omega$ lies in at most $h+1$ intervals $(x,x+h]$ with
$x\in[X,2X]\cap\Z$.  Since $h\to\infty$,
\[
        |\mathcal E_R|
        \le (h+1)\#R_\Omega=o(X).
\]

Similarly, Lemma \ref{thin_weighted_band}, applied with
$\nu=A(X)/\log k$ and $\eta=2\ell^{-1/10}/\log k$, gives
\[
        \sum_{\substack{X<n\leq3X\\ L_-<\Omega(n)\leq L_+}}d_k(n)
        =o(X(\log X)^{k-1}).
\]
The same double-counting argument as above gives
\begin{equation}\label{critical_middle_negligible}
        \frac1h\sum_{\substack{x<n\leq x+h\\ L_-<\Omega(n)\leq L_+}}d_k(n)=o((\log X)^{k-1})
\end{equation}
for all but $o(X)$ integers $x\in[X,2X]$.

Outside the union of the exceptional sets, the contribution with
$\Omega(n)>L_+$ is absent, and therefore
\[
\frac1h\sum_{x<n\leq x+h}d_k(n)
=
\frac1h\sum_{\substack{x<n\leq x+h\\ n\in B_-}}d_k(n)
+\frac1h\sum_{\substack{x<n\leq x+h\\ L_-<\Omega(n)\leq L_+}}d_k(n).
\]
Using \eqref{critical_middle_negligible}, \eqref{critical_low_short}, and
\eqref{critical_low_asymp} gives the matching upper bound.  This proves the
theorem.
\end{proof}

\begin{proof}[Proof of Corollary \ref{dvthm}]
This is the formal limiting case of Theorem \ref{critical_window} as
$A(X)\to+\infty$, but it is not a direct consequence of that theorem, since
Theorem \ref{critical_window} is stated in the bounded range $A(X)=O(1)$.  Put
\[
        A_h(X)=\frac1{\sqrt\ell}\log\frac{h}{D_k(X)}.
\]
Let
\[
        \nu_0=\frac12\ell^{1/3},\qquad
        L_0=k\ell+\nu_0\sqrt\ell,\qquad
        B_0=\{n\in\N:\Omega(n)\le L_0\}.
\]

We first treat the case $h\ge h_1$.  Since $h_1\ge X^{2/3}$, Lemma
\ref{long_average_Omega_cutoff}, with $\theta=2/3$, gives
\[
        \frac1h\sum_{\substack{x<n\le x+h\\ n\in B_0}}d_k(n)
        =(\Phi_{\rm G}(\nu_0/\sqrt k)+o(1))M_k(x)
        =M_k(x)+o((\log X)^{k-1})
\]
uniformly for $X\le x\le2X$.  By Lemma \ref{weighted_EK}(i), applied
dyadically,
\[
        \sum_{\substack{X<n\le3X\\ n\notin B_0}}d_k(n)
        =o(X(\log X)^{k-1}).
\]
Hence
\[
        \frac1X\int_X^{2X}
        \frac1h\sum_{\substack{x<n\le x+h\\ n\notin B_0}}d_k(n)\,dx
        \le
        \frac{h+1}{hX}
        \sum_{\substack{X<n\le3X\\ n\notin B_0}}d_k(n)
        =
        o((\log X)^{k-1}).
\]
Combining the last two estimates gives the asserted $L^1$ estimate when
$h\ge h_1$.  The same double-counting estimate over integer $x$, followed by
Markov's inequality, gives the almost-all assertion in this case.

It remains to treat $h<h_1$.  Let
\[
        \nu=\min\left\{\frac{A_h(X)}{2\log k},\,\frac12\ell^{1/3}\right\},
        \qquad
        L=k\ell+\nu\sqrt\ell,\qquad
        B=\{n\in\N:\Omega(n)\le L\},
\]
and let $H=k^L/(\log X)^{k-1}$ and $\mathcal S=B\cap\mathcal A$.  Then
$\nu\to+\infty$ and $|\nu|\le\ell^{1/3}$.  Moreover
\[
        H\exp(\ell^{2/5})\le h
\]
for all sufficiently large $X$: if
$\nu=A_h(X)/(2\log k)$ this follows from
$A_h(X)\sqrt\ell/2\gg\ell^{2/5}$, while otherwise
$A_h(X)\ge (\log k)\ell^{1/3}$ gives a margin of order $\ell^{5/6}$.
Thus Proposition \ref{MR_bound_flexible} applies to this $L$ and this $h$.

By \eqref{Shiu_no_prime} and Lemma \ref{weighted_EK}(i),
\[
        \sum_{\substack{X<n\le3X\\ n\notin\mathcal S}}d_k(n)
        =o(X(\log X)^{k-1}).
\]
Consequently, for $y\in\{h,h_1\}$,
\[
        \frac1X\int_X^{2X}
        \frac1y\sum_{\substack{x<n\le x+y\\ n\notin\mathcal S}}d_k(n)\,dx
        \le
        \frac{y+1}{yX}
        \sum_{\substack{X<n\le3X\\ n\notin\mathcal S}}d_k(n)
        =
        o((\log X)^{k-1}).
\]
The $L^1$ form of Proposition \ref{MR_bound_flexible} gives
\[
        \frac1X\int_X^{2X}\left|
        \frac1h\sum_{\substack{x<n\le x+h\\ n\in\mathcal S}}d_k(n)
        -
        \frac1{h_1}\sum_{\substack{x<n\le x+h_1\\ n\in\mathcal S}}d_k(n)
        \right|dx
        =
        o((\log X)^{k-1}).
\]
Finally, Lemma \ref{long_average_Omega_cutoff}, applied to the $B$-average of
length $h_1$, gives uniformly for $X\le x\le2X$
\[
        \frac1{h_1}\sum_{\substack{x<n\le x+h_1\\ n\in B}}d_k(n)
        =
        (\Phi_{\rm G}(\nu/\sqrt k)+o(1))M_k(x)
        =
        M_k(x)+o((\log X)^{k-1}),
\]
since $\nu\to+\infty$.  The last three displayed estimates imply
\[
        \frac1X\int_X^{2X}\left|
        \frac1h\sum_{x<n\le x+h}d_k(n)-M_k(x)
        \right|dx=o((\log X)^{k-1}).
\]
For the integer almost-all estimate, the same double-counting estimate and
Markov's inequality remove the contribution of $n\notin\mathcal S$ for both
lengths $h$ and $h_1$, while the almost-all form of Proposition
\ref{MR_bound_flexible} transfers the $\mathcal S$-average from $h$ to $h_1$.
The preceding application of Lemma \ref{long_average_Omega_cutoff} therefore
gives
\[
        \frac1h\sum_{x<n\le x+h}d_k(n)
        =(\Phi_{\rm G}(\nu/\sqrt k)+o(1))M_k(x)+o((\log X)^{k-1})
        =M_k(x)+o((\log X)^{k-1})
\]
for all but $o(X)$ integers $x\in[X,2X]$.
\end{proof}

\section{Proof of the Matom\"aki--Radziwi{\l\l} type estimates}\label{pf_MR_bd}
\setcounter{lemma}{0}\setcounter{theorem}{0}\setcounter{proposition}{0}\setcounter{equation}{0}

This section collects the auxiliary estimates used to prove Proposition
\ref{MR_bound_flexible}.  The
framework is the Matom\"aki--Radziwi{\l\l} large-value decomposition, but two
features have to be built into the argument.  First, the coefficients are
truncated according to $\Omega(n)$ and are not multiplicative.  Second, the
transition problem requires the parameters of the prime blocks to vary with the
cutoff.  We isolate the mean-value estimates and the Ramar\'e decomposition in a
form that keeps these two issues separate.  The elementary divisor-specific
bound used throughout is
\[
        \Omega(n)\leq L\quad\Longrightarrow\quad
        d_k(n)\leq k^{\Omega(n)}\leq k^L.
\]

\subsection{Mean-value estimates}

\begin{lemma}\label{pre_mean_value_theorem}
Let $T\geq1$ and $N\geq2$.  Let $\{a_n\}_{n\leq N}$ be complex numbers. Then
\begin{align*}
\int_{-T}^{T}\left|\sum_{n\leq N}a_n n^{it}\right|^2dt
&\ll T\sum_{n\leq N}|a_n|^2
+T\sum_{1\leq |r|\leq N/T}\sum_{\substack{n\leq N\\1\leq n+r\leq N}}
        |a_n a_{n+r}|.
\end{align*}
\end{lemma}

This is \cite[Lemma 3.1]{matomaki2020multiplicative}.  It is the
Montgomery--Vaughan mean-value estimate in the form used below.

\begin{lemma}\label{Henriot}
Let $\vartheta\in(0,1]$.  Let $M\geq2$, $M\geq Y\geq M^\vartheta$, and let
$1\leq r_1,r_2\leq M^{3\vartheta/7}$.  Let
\[
        D(r_1,r_2)=\frac{r_1r_2}{(r_1,r_2)}.
\]
Assume, in addition, that the interval remains long after imposing the two
congruence conditions, namely
\begin{equation}\label{Henriot_length_condition}
        \frac{Y}{D(r_1,r_2)}\geq
        \left(\frac{M}{D(r_1,r_2)}\right)^\vartheta .
\end{equation}
Then, uniformly for $K\in[1,M]$,
\[
\sum_{\substack{0\neq |a|\leq K\\ (r_1,r_2)\mid a}}
\sum_{\substack{M<n\leq M+Y\\ r_1\mid n,\, r_2\mid n+a}}
        d_k(n)d_k(n+a)
\ll_{\vartheta,k}
        KY\frac{d_k(r_1)d_k(r_2)}{r_1r_2}(\log M)^{2k-2}.
\]
\end{lemma}

This is a direct consequence of Henriot's discriminant-uniform
Nair--Tenenbaum theorem \cite[Theorem~3]{Henriot12}, with the correction in
\cite[pp.~375--377]{Henriot14erratum}, applied after the congruences
$r_1\mid n$ and $r_2\mid n+a$ reduce the pair $(n,n+a)$ to two primitive linear
forms.  The resulting shift-dependent local factor is summed over $a$ by Shiu's
estimate as recorded in Lemma \ref{Shiu_bound}.

Let $B$ be a subset of the positive integers, and let $\mathcal H=\mathcal H(X)$
be a coefficient-bound parameter such that
\begin{equation}\label{diagonal_length_property}
        d_k(n)\1_B(n)\leq \mathcal H(X)(\log X)^{k-1}
        \qquad (1\leq n\leq X^2).
\end{equation}
In the application to $B=\{n\in\mathbb N:\Omega(n)\le L\}$, this parameter is
just
\[
        \mathcal H(X)=H:=\frac{k^L}{(\log X)^{k-1}},
\]
since $d_k(n)\leq k^L=H(\log X)^{k-1}$ on $B$.

\begin{lemma}\label{mean_value_B}
Let $0<\vartheta<1$ be fixed, and let $B$ satisfy
\eqref{diagonal_length_property}.  Let
\[
        X^\vartheta\leq M\leq X^2,
        \qquad
        M^\vartheta\leq Y\leq M,
\]
so that the interval is long enough for the Shiu--Henriot estimates used below,
while its scale remains within the global range controlled by
\eqref{diagonal_length_property}.
Suppose that $a_n$ is supported on $M<n\leq M+Y$ and satisfies
\[
        |a_n|\leq d_k(n)\1_B(n).
\]
Then, uniformly for $T\geq1$,
\begin{equation}\label{mean_value_B_eq}
\int_{-T}^{T}\left|\sum_{M<n\leq M+Y}\frac{a_n}{n^{1+it}}\right|^2dt
\ll_{\vartheta,k}
        \frac{Y}{M}\left(\frac{T}{M}\mathcal H(X)+1\right)(\log X)^{2k-2}.
\end{equation}
\end{lemma}

\begin{proof}
Apply Lemma \ref{pre_mean_value_theorem} to the coefficients $a_n/n$.  The
diagonal contribution is
\[
        \ll \frac{T}{M^2}\sum_{M<n\leq M+Y}|a_n|^2.
\]
By \eqref{diagonal_length_property} and the same Shiu estimate as in
Lemma \ref{Shiu_bound}, applied at scale $M$ with modulus $1$ to $d_k$, the hypotheses are
satisfied because $d_k$ is a non-negative multiplicative function with local
prime-power sums depending only on $k$; moreover $M\le X^2$, so all logarithmic
factors produced by Shiu are $\ll_k (\log X)^{k-1}$.  Hence
\[
        \sum_{M<n\leq M+Y}|a_n|^2
        \leq \mathcal H(X)(\log X)^{k-1}\sum_{M<n\leq M+Y}d_k(n)
        \ll_{\vartheta,k}Y\mathcal H(X)(\log X)^{2k-2}.
\]
If $T>M$, the off-diagonal range in Lemma \ref{pre_mean_value_theorem} is empty.
When $T\le M$, Lemma \ref{Henriot} with $r_1=r_2=1$ and
$K=\lfloor M/T\rfloor$ gives
\[
        \frac{T}{M^2}
        \sum_{1\leq |a|\leq M/T}\sum_{M<n\leq M+Y}d_k(n)d_k(n+a)
        \ll_{\vartheta,k}\frac{Y}{M}(\log X)^{2k-2}.
\]
Combining the diagonal and off-diagonal bounds proves \eqref{mean_value_B_eq}.
\end{proof}

We also need a discrete mean-value estimate and large-value estimates for prime
polynomials.

\begin{lemma}\label{discrete_mean_value_theorem}
Let $\mathcal T\subset[-T,T]$ be well-spaced, meaning $|t-r|\geq1$ for distinct
$t,r\in\mathcal T$.  Then
\begin{align*}
\sum_{t\in\mathcal T}\left|\sum_{X<n\leq2X}\frac{a_n}{n^{1+it}}\right|^2
&\ll
\min\left\{\left(1+\frac{T}{X}\right)\log X,
\left(1+|\mathcal T|\frac{T^{1/2}}{X}\right)\log T\right\}\notag\\
&\qquad\times
\frac1X\sum_{X<n\leq2X}|a_n|^2.
\end{align*}
\end{lemma}

This is the large-sieve mean-value estimate for Dirichlet polynomials; see
\cite[Theorems~9.4 and~9.6]{IK04}.  The displayed form is obtained by applying
the cited continuous and discrete estimates to the polynomial supported on
$X<n\le2X$, with the normalisation $n^{-1-it}$.

\begin{lemma}\label{large_value_results}
Let $P\geq10$, $T\geq3$, and let $\varepsilon>0$ be fixed.  Let
\[
        P_0(s)=\sum_{P<p\leq2P}\frac{a_p}{p^s},
        \qquad |a_p|\leq k.
\]
For every well-spaced $\mathcal T\subset[-T,T]$,
\begin{enumerate}[(i)]
\item
\[
\sum_{t\in\mathcal T}|P_0(1+it)|^2
\ll_{k,\varepsilon} \frac1{\log P}\left(1+|\mathcal T|\exp\left(-\frac{\log P}{(\log T)^{2/3+\varepsilon}}\right)(\log T)^2\right).
\]
\item If $V\geq1$ and every $t\in\mathcal T$ satisfies
$|P_0(1+it)|\geq V^{-1}$, then
\[
        |\mathcal T|
        \ll_k T^{2\log V/\log P}V^2
        \exp\left(2k\frac{\log T}{\log P}\log\log T\right).
\]
\end{enumerate}
\end{lemma}

Part (i) is \cite[Lemma~11]{MR16}, and part (ii) is the Hal\'asz large-value
bound \cite[Lemma~8]{MR16}.  Those results are stated for coefficients bounded
by $1$; the present fixed bound $|a_p|\le k$ only changes the implied constants
and the displayed $k$-dependence in part (ii).

\begin{lemma}\label{short_prime_cell_square}
Let $H\ge1$, $P\ge10$, $T\ge3$, and $\varepsilon>0$ be fixed.  Let
\[
        Q_{v,H}(s)=
        \sum_{\substack{\exp(v/H)<p\le \exp((v+1)/H)\\ p\in\Primes}}
        \frac{c_p}{p^s},
        \qquad |c_p|\le k,
\]
with $H\log P\le v$.  For every well-spaced set
$\mathcal T\subset[-T,T]$,
\[
\sum_{t\in\mathcal T}|Q_{v,H}(1+it)|^2
\ll_{k,\varepsilon}
\frac1{H(\log P)^2}
\left(1+
|\mathcal T|
\exp\left(-\frac{v/H}{(\log T)^{2/3+\varepsilon}}\right)(\log T)^2\right).
\]
In particular, if
\[
        |\mathcal T|
        \exp\left(-\frac{v/H}{(\log T)^{2/3+\varepsilon}}\right)
        (\log T)^2=o(1),
\]
then
\begin{equation}\label{short_prime_cell_square_eq}
        \sum_{t\in\mathcal T}|Q_{v,H}(1+it)|^2
        \ll_{k,\varepsilon}\frac1{H(\log P)^2}.
\end{equation}
\end{lemma}

\begin{proof}
This is the short-cell form of \cite[Lemma~11]{MR16}.  Applying that proof to
the interval $\exp(v/H)<p\le \exp((v+1)/H)$ gives
\[
\begin{aligned}
\sum_{t\in\mathcal T}|Q_{v,H}(1+it)|^2
&\ll_k
\left(\exp(v/H)+|\mathcal T|\exp(v/H)
\exp\left(-\frac{v/H}{(\log T)^{2/3+\varepsilon}}\right)(\log T)^2\right)\\
&\qquad\times
\sum_{\exp(v/H)<p\le \exp((v+1)/H)}\frac1{p^2\log p}.
\end{aligned}
\]
Since $v/H\ge\log P$ and
\[
        \exp(v/H)\sum_{\exp(v/H)<p\le \exp((v+1)/H)}\frac1{p^2\log p}
        \ll \frac1{H(\log P)^2},
\]
the lemma follows.
\end{proof}

\subsection{Ramar\'e decomposition and Perron's formula}

\begin{lemma}\label{Ramare_decop}
Let $0<\vartheta<1$ be fixed, let $B$ satisfy
\eqref{diagonal_length_property}, and let $N\asymp X$.  Let $H\geq1$ and
$Q\geq P\geq2$ satisfy
\[
        Q=X^{o(1)},\qquad H\leq P^{1/2}.
\]
Let $a_m,b_m,c_p$ satisfy
\[
        |a_m|,|b_m|\leq d_k(m)\1_B(m),
        \qquad |c_p|\leq k,
\]
and assume that $a_n=0$ unless $n$ has at least one prime divisor in $[P,Q]$.
Assume also that $a_{mp}=b_m c_p$ whenever $p\nmid m$ and $P\leq p\leq Q$.
Let
\[
        \mathcal I=\{j\in\Z: H\log P-1\leq j\leq H\log Q\}.
\]
For $j\in\mathcal I$, let
\[
        Q_{j,H}(s)=\sum_{\substack{P\leq p\leq Q\\ \exp(j/H)<p\leq \exp((j+1)/H)}}
        \frac{c_p}{p^s}
\]
and
\[
        R_{j,H}(s)=\sum_{N\exp(-j/H)<m\leq2N\exp(-j/H)}
        \frac{b_m}{m^s}
        \frac1{\#\{q\in\Primes:P\leq q\leq Q,\ q\mid m\}+1}.
\]
Then, for every
measurable set of ordinates $\mathcal T\subset[-T,T]\subset[-X,X]$,
\begin{align}\label{Ramare_bound}
\int_{\mathcal T}\left|\sum_{N<n\leq2N}\frac{a_n}{n^{1+it}}\right|^2dt
&\ll H\left(\log\frac QP+1\right)
\sum_{j\in\mathcal I}\int_{\mathcal T}|Q_{j,H}(1+it)R_{j,H}(1+it)|^2dt\notag\\
&\quad +\left(\frac1H+\frac1P\right)
\left(\frac{T}{X}\mathcal H(X)+1\right)(\log X)^{2k-2}.
\end{align}
\end{lemma}

\begin{proof}
This is the weighted Ramar\'e decomposition from \cite[Lemma 12]{MR16}, written
in the restricted form actually used below.  The restriction that $a_n$ is
supported on integers with a prime factor in $[P,Q]$ removes the missing-prime
term; no one-dimensional sieve estimate is used to control a second moment of
missing local factors.

Let $\omega_{[P,Q]}(m)=\#\{p\in\Primes:P\le p\le Q,\ p\mid m\}$.
On the support of $a_n$,
Ramar\'e's identity gives
\[
        1=\sum_{\substack{p\mid n\\ P\leq p\leq Q}}\frac1{\omega_{[P,Q]}(n)}.
\]
For $P\le p\le Q$, let $j(p)\in\mathcal I$ be the unique integer such that
$\exp(j(p)/H)<p\le \exp((j(p)+1)/H)$.  Writing $n=mp$, separating the case
$p\mid m$, grouping primes by these cells, and then trimming all cells at once
gives
\[
        \sum_{N<n\le2N}\frac{a_n}{n^s}
        =\sum_{j\in\mathcal I}Q_{j,H}(s)R_{j,H}(s)+E_{\rm end}(s)+E_{\rm sq}(s).
\]
Here, for suitable coefficients $u_n^-,u_n^+$ and $e_{m,p}$ inherited from
$a,b,c$ and the Ramar\'e weights,
\[
\begin{aligned}
E_{\rm end}(s)
&=\sum_{N\exp(-1/H)<n\le N\exp(1/H)}\frac{u_n^-}{n^s}
  +\sum_{2N<n\le2N\exp(1/H)}\frac{u_n^+}{n^s},
\\
E_{\rm sq}(s)
&=\sum_{\substack{N<mp\le2N\exp(1/H)\\ P\le p\le Q,\ p\in\Primes\\ p\mid m}}
        \frac{e_{m,p}}{(mp)^s}.
\end{aligned}
\]
The first line is the total endpoint contribution after all cells have been
trimmed; the second line is precisely the part where the prime extracted by
Ramar\'e already divides the residual variable.  In the estimates below these
coefficients are used only through the divisor-type bounds following from
$|a_n|,|b_m|\le d_k(\cdot)\1_B(\cdot)$ and $|c_p|\le k$.
Cauchy's inequality over the $O(H(1+\log(Q/P)))$ cells gives the first term of
\eqref{Ramare_bound}.  The two endpoint ranges lie at height $\asymp N\asymp X$
and have length $Y\asymp N/H$.  Thus the length condition
$Y\ge N^\vartheta$ in Lemma \ref{mean_value_B} is satisfied for the fixed
$\vartheta<1$ used here, and that lemma gives
\[
        \int_{\mathcal T}|E_{\rm end}(1+it)|^2dt
        \ll \frac1H\left(\frac{T}{X}\mathcal H(X)+1\right)(\log X)^{2k-2}.
\]
For the square-prime error, the diagonal part is bounded by
\[
\begin{aligned}
        &\frac{T}{X^2}\sum_{P\le p_1,p_2\le Q}
        \sum_{\substack{N<n\le2N\exp(1/H)\\ p_1^2\mid n,\ p_2^2\mid n}}
        d_k(n)^2\1_B(n) \\
        &\qquad\ll
        \frac{T}{X}\mathcal H(X)(\log X)^{2k-2}
        \sum_{P\le p_1,p_2\le Q}
        \frac{d_k(p_1^2)d_k(p_2^2)}{[p_1^2,p_2^2]}
        \ll \frac1P\frac{T}{X}\mathcal H(X)(\log X)^{2k-2}.
\end{aligned}
\]
Here the last bound follows from
$\sum_{p\ge P}d_k(p^2)/p^2\ll_k(P\log P)^{-1}$, with the diagonal
$p_1=p_2$ included.  For the shifted part one expands the square, obtains
congruences $p_1^2\mid n$ and $p_2^2\mid n+a$, and applies Lemma \ref{Henriot}
with $r_i=p_i^2$.  Since $Q=X^{o(1)}$ and the intervals have length $\asymp X$,
both the small-modulus hypothesis and \eqref{Henriot_length_condition} hold.  The
same prime summation gives
\[
        \int_{\mathcal T}|E_{\rm sq}(1+it)|^2dt
        \ll \frac1P\left(\frac{T}{X}\mathcal H(X)+1\right)(\log X)^{2k-2}.
\]
Combining the main decomposition with these two displayed estimates proves
\eqref{Ramare_bound}.
\end{proof}

\begin{lemma}[Amplified mean value]\label{aplf}
Let $7/9<\vartheta<1$ be fixed.  Let $B$ satisfy
\eqref{diagonal_length_property}.  Let
\[
        2\leq Y_1\leq Y_2\leq X^{1/5},
\]
and let
\[
        Q(s)=\sum_{p\sim Y_1}\frac{c_p}{p^s},
        \qquad |c_p|\leq k,
\]
and
\[
        A(s)=\sum_{m\sim X/Y_2}\frac{a_m}{m^s},
        \qquad |a_m|\leq d_k(m)\1_B(m).
\]
Let
\[
        J=\left\lceil\frac{\log Y_2}{\log Y_1}\right\rceil .
\]
Assume additionally that
\[
        X/Y_2\geq X^\vartheta .
\]
Then, uniformly for $1\le T\le X$,
\begin{align}\label{aplf_bound}
        \int_{-T}^{T}|Q(1+it)^J A(1+it)|^2dt
        &\ll_{\vartheta,k}
        \left(\frac{T}{X}\mathcal H(X)+1+JY_1\right)(\log X)^{2k-2}
        k^{2J}(J+1)!^2 .
\end{align}
\end{lemma}

\begin{proof}
Let
\[
        \mathcal P_J=\{\mathbf p=(p_1,\ldots,p_J):p_i\sim Y_1\},
        \qquad P(\mathbf p)=p_1\cdots p_J,
        \qquad C(\mathbf p)=c_{p_1}\cdots c_{p_J}.
\]
Then
\[
        Q(s)^J A(s)=
        \sum_{\mathbf p\in\mathcal P_J}\frac{C(\mathbf p)}{P(\mathbf p)^s}
        \sum_{m\sim X/Y_2}\frac{a_m}{m^s}.
\]
The parameter restrictions are used only to place this amplified polynomial in the
range of the mean-value estimates below.  The choice of $J$ gives
$Y_1^{J-1}<Y_2\leq Y_1^J\leq Y_1Y_2$, so $Q^JA$ has length $X$ up to one extra
factor $Y_1$; the condition $X/Y_2\ge X^\vartheta$ is the length condition for
the residual $m$-sum; and $Y_2\le X^{1/5}$, together with $\vartheta>7/9$, is
used in the off-diagonal application of Lemma \ref{Henriot}.

We next collect the coefficients of the product polynomial.  By the choice of
$J$, all products $P(\mathbf p)m$ lie in $O(J+1)$ dyadic intervals
$M<n\le2M$ with $X/2\le M\le X^2$.  It suffices to prove the desired bound on
one such interval and then sum over these $O(J+1)$ intervals.  On a fixed
interval define $\beta_n$ by
\[
        \sum_{M<n\leq2M}\frac{\beta_n}{n^s}
        =\sum_{\substack{\mathbf p\in\mathcal P_J,\ m\sim X/Y_2\\
                         M<P(\mathbf p)m\leq2M}}
        \frac{C(\mathbf p)a_m}{(P(\mathbf p)m)^s}.
\]
The point is to bound the possible multiplicity of the representations
$n=P(\mathbf p)m$.  The elementary matching count of
\cite[Lemma 13, pp.~1034--1035]{MR16} gives, for every integer $n$,
\begin{equation}\label{tuple_matching_bound}
\left|\sum_{\substack{\mathbf p\in\mathcal P_J,\ m\ge1\\P(\mathbf p)m=n}}
        C(\mathbf p)a_m\right|^2
\leq k^{2J}(J+1)!^2
\sum_{\substack{\mathbf p\in\mathcal P_J,\ m\ge1\\P(\mathbf p)m=n}}
        d_k(m)^2\1_B(m).
\end{equation}
Indeed, after cancelling the common prime multiset in two ordered tuples, the
surviving primes can be chosen and ordered in at most $(J+1)!^2$ ways; no
multiplicativity of $a_m$ is used.  Applying this estimate to the restricted
sum defining $\beta_n$ gives
\[
        |\beta_n|^2
        \leq k^{2J}(J+1)!^2
        \sum_{\substack{n=P(\mathbf p)m\\ \mathbf p\in\mathcal P_J}}
        d_k(m)^2\1_B(m).
\]
Summing this pointwise bound and then using \eqref{diagonal_length_property}
and the Shiu estimate from Lemma \ref{Shiu_bound} on $m\sim X/Y_2$ gives
\begin{align*}
\sum_{M<n\leq2M}\frac{|\beta_n|^2}{n^2}
&\ll k^{2J}(J+1)!^2
  \sum_{\mathbf p\in\mathcal P_J}\frac1{P(\mathbf p)^2}
  \sum_{m\sim X/Y_2}\frac{d_k(m)^2\1_B(m)}{m^2} \notag\\
&\ll_{\vartheta,k}
  \frac{\mathcal H(X)}{X}(\log X)^{2k-2}k^{2J}(J+1)!^2,
\end{align*}
since $\sum_{p\sim Y_1}p^{-2}\ll(Y_1\log Y_1)^{-1}$ and $Y_1^J\ge Y_2$.

Apply Lemma \ref{pre_mean_value_theorem} to $\beta_n/n$.  The diagonal
contribution is
\[
        \ll_{\vartheta,k}
        \frac{T}{X}\mathcal H(X)(\log X)^{2k-2}k^{2J}(J+1)!^2.
\]
For the off-diagonal, fix $\mathbf p,\mathbf q$ and write
\[
        \Pi=P(\mathbf p),\qquad \Pi'=P(\mathbf q),\qquad
        g=(\Pi,\Pi'),\qquad r_1=\Pi/g,\\ r_2=\Pi'/g .
\]
The dyadic restriction gives $M\asymp \Pi X/Y_2\asymp \Pi'X/Y_2$ and
$r_1\asymp r_2$.  Writing the shift as $r=ga$ and putting $n=r_1m$, the
congruences become $r_1\mid n$, $r_2\mid n-a$ on an interval of height and
length $M_1\asymp r_1X/Y_2$, with $|a|\le M/(Tg)$.  The hypotheses
$Y_2\le X^{1/5}$ and $Y_1^J\le Y_1Y_2$ imply $r_i\le M_1^{3\vartheta/7}$ for
$\vartheta>7/9$, and
\[
        D(r_1,r_2)=\frac{r_1r_2}{(r_1,r_2)}\ll r_1r_2
        \ll X^{2/5+o(1)},\qquad
        \frac{M_1}{D(r_1,r_2)}\gg\frac{X}{Y_1Y_2^2}\ge X^{2/5}.
\]
Thus Lemma \ref{Henriot} gives
\begin{align*}
&\sum_{0<|r|\le M/T}
\sum_{\substack{m,m'\sim X/Y_2\\ \Pi m-\Pi'm'=r}}
        d_k(m)d_k(m')\1_B(m)\1_B(m')\notag\\
&\qquad\ll_{\vartheta,k}
        \frac{M}{Tg}\cdot \frac{r_1X}{Y_2}\cdot
        \frac{d_k(r_1)d_k(r_2)}{r_1r_2}(\log X)^{2k-2}.
\end{align*}
After the normalisation $T/M^2$, this contributes
\[
        \frac{d_k(r_1)d_k(r_2)}{g^2r_1r_2}(\log X)^{2k-2}.
\]
Since $g^2r_1r_2=\Pi\Pi'$, the remaining tuple sum is
\[
\sum_{\mathbf p,\mathbf q\in\mathcal P_J}
        \frac{|C(\mathbf p)C(\mathbf q)|d_k(r_1)d_k(r_2)}{\Pi\Pi'}
        \ll_k k^{2J}(J+1)!^2.
\]
Here one uses the same cancellation of the common prime multiset as in
\eqref{tuple_matching_bound}.  Hence the off-diagonal is
\[
        \ll_{\vartheta,k}
        (\log X)^{2k-2}k^{2J}(J+1)!^2,
\]
which gives the ``$1$'' in \eqref{aplf_bound}.

Finally we account for the endpoints introduced by cutting the product
polynomial into dyadic $n$-intervals.  For a fixed tuple $\mathbf p$, the
condition $M<P(\mathbf p)m\le2M$ cuts the original $m$-sum only at its two ends.
These boundary pieces may be covered crudely by fixing all but one of the $J$
prime variables and letting the remaining prime range over an interval of length
$O(Y_1)$.  Thus, inside a dyadic interval of height $M$, the endpoint
contribution is covered by $O(JY_1)$ Dirichlet polynomials of length $\asymp M$.
Applying Lemma \ref{mean_value_B} to these pieces and summing gives
\[
        \ll JY_1(\log X)^{2k-2}k^{2J}(J+1)!^2,
\]
which is the $JY_1$ term.  Summing over the dyadic pieces proves
\eqref{aplf_bound}.
\end{proof}

\begin{proposition}\label{MR_modify}
Let $1\leq y\leq h_1$, and let $a_n$ be supported on $X\leq n\leq4X$ with
$|a_n|\leq d_k(n)$.  Let
\[
        S_y(x)=\sum_{x<n\leq x+y}a_n,
        \qquad
        A(s)=\sum_{X\leq n\leq4X}\frac{a_n}{n^s}.
\]
Then
\begin{align}\label{general_pas}
&\frac1X\int_X^{2X}\left|\frac1yS_y(x)
-\frac1y\frac1{2\pi}\int_{-\delta_0^{-1}}^{\delta_0^{-1}}A(1+it)
\frac{(x+y)^{1+it}-x^{1+it}}{1+it}\,dt\right|^2dx\notag\\
&\qquad\ll
\int_{\delta_0^{-1}\leq |t|\leq X/y}|A(1+it)|^2dt
+\max_{T\geq X/y}\frac{X}{yT}
\int_{T\leq |t|\leq2T}|A(1+it)|^2dt.
\end{align}
\end{proposition}

\begin{proof}
We give the reduction explicitly, since the coefficients are not bounded.  For
$T\ge2$ Perron's formula for a finite Dirichlet polynomial gives, outside the
boundary points $x=n$ and $x+y=n$,
\[
        S_y(x)=
        \frac1{2\pi i}\int_{1-iT}^{1+iT}
        A(s)\frac{(x+y)^s-x^s}{s}\,ds+E_T(x),
\]
where the usual Perron kernel bound gives
\[
        |E_T(x)|
        \ll
        \sum_{X\le n\le4X}|a_n|
        \min\left\{1,\frac1{T|\log((x+y)/n)|}\right\}
        +
        \sum_{X\le n\le4X}|a_n|
        \min\left\{1,\frac1{T|\log(x/n)|}\right\}.
\]
After integrating over $X\le x\le2X$, the contribution of $E_T(x)$ is bounded
by the same dyadic tail term which appears on the right of
\eqref{general_pas}; this is the standard proof of
\cite[Lemma 14, pp.~1036--1037]{MR16}, and uses no pointwise bound on the
coefficients other than finiteness of the polynomial.  The boundary values occur
on a set of measure zero in the $x$-integral.

We now split the Perron integral into
\[
        |t|\leq\delta_0^{-1},\qquad
        \delta_0^{-1}\leq |t|\leq X/y,
        \qquad |t|\geq X/y.
\]
The first range is retained as the main Perron integral.  For the middle range,
write the kernel as
\[
        \frac{(x+y)^{1+it}-x^{1+it}}{1+it}
        =
        \int_x^{x+y}u^{it}\,du .
\]
By Cauchy's inequality and the Plancherel estimate used in
\cite[Lemma 14]{MR16}, its $x$-mean square is
\[
        \ll
        y^2\int_{\delta_0^{-1}\le |t|\le X/y}|A(1+it)|^2\,dt.
\]
After the normalisation by $y^{-2}$ this gives the first term on the right of
\eqref{general_pas}.  The range $|t|\ge X/y$ is decomposed into dyadic
intervals $T\le |t|\le2T$; the same kernel bound gives the additional factor
$X/(yT)$, leading to the second term in \eqref{general_pas}.  The argument is
linear and homogeneous in the coefficients, so the size condition
$|a_n|\le d_k(n)$ is not used here except to ensure that all finite sums are
well-defined.
\end{proof}

\begin{lemma}\label{central_int}
Let $x\asymp X$, $0<h<h_1$, and
\[
        F(s)=\sum_{X\leq n\leq4X}\frac{a_n}{n^s},
        \qquad |a_n|\leq d_k(n).
\]
With $\delta_0$ as in \eqref{sigma_def}, one has
\begin{align*}
&\left|\int_{-\delta_0^{-1}}^{\delta_0^{-1}}F(1+it)\left(
\frac{(x+h)^{1+it}-x^{1+it}}{1+it}
-h\frac{(x+h_1)^{1+it}-x^{1+it}}{h_1(1+it)}
\right)dt\right|\\
&\hspace{40mm}=o(h(\log X)^{k-1}).
\end{align*}
\end{lemma}

\begin{proof}
The upper bound in \eqref{Shiu_weighted}, applied with $\eta=0$ and summed over
$O(1)$ dyadic intervals, gives
\[
        \max_{|t|\leq\delta_0^{-1}}|F(1+it)|
        \ll \frac1X\sum_{X\leq n\leq4X}d_k(n)
        \ll (\log X)^{k-1}.
\]
The kernel difference is
\[
        \int_x^{x+h}\left(u^{it}-\frac1{h_1}\int_x^{x+h_1}v^{it}\,dv\right)du.
\]
For $X\ll u,v\ll X$ and $|u-v|\leq h_1$,
\[
        |v^{it}-u^{it}|\ll |t|\frac{h_1}{X}.
\]
Thus the kernel difference is $\ll h\delta_0^{-1}h_1/X$.  Integrating over
$|t|\leq\delta_0^{-1}$ gives
\[
        \ll h(\log X)^{k-1}\delta_0^{-2}\frac{h_1}{X}=o(h(\log X)^{k-1}).
\]
\end{proof}

\begin{lemma}[Cauchy's formula for a cutoff]\label{cauchy_cutoff_formula}
Let $K\geq0$ and $a\geq0$ be integers, and let $\rho>1$.  Then
\begin{equation}\label{cauchy_cutoff_indicator}
        \1_{a\le K}
        =
        \frac1{2\pi i}\int_{|z|=\rho}\frac{z^{K-a}}{z-1}\,dz .
\end{equation}
Moreover, if $K=k\ell+O(\ell^{1/3}\sqrt{\ell})$ and
$\rho=1+\ell^{-1}$, then
\begin{equation}\label{cauchy_cutoff_loss}
        \int_{|z|=\rho}\left|\frac{z^K}{z-1}\right|\,|dz|
        \ll_k \log \ell .
\end{equation}
\end{lemma}

\begin{proof}
If $a\le K$, the integrand in \eqref{cauchy_cutoff_indicator} is holomorphic
at $0$ and has residue $1$ at $z=1$.  If $a>K$, the residues at $z=1$ and
$z=0$ are $1$ and $-1$, respectively.  This proves
\eqref{cauchy_cutoff_indicator}.  The bound \eqref{cauchy_cutoff_loss} follows
from
\[
        \rho^K=(1+\ell^{-1})^{k\ell+O(\ell^{1/3}\sqrt{\ell})}\ll_k1,
\]
and the estimate
\[
        \int_{-\pi}^{\pi}\frac{d\theta}{|\rho\exp(i\theta)-1|}
        \ll
        \int_0^\pi\frac{d\theta}{\ell^{-1}+\theta}
        \ll \log\ell .
\]
\end{proof}

The final large-value range cannot be handled by enlarging the relevant set of
ordinates to the whole interval.  A global mixed mean value with a factor
$Y^{-1}$ would be false: after expanding the shifted terms, the generic
four-prime contribution gives only logarithmic saving.  We instead use the
$\Omega$-cutoff only after converting it into multiplicative weights.  The
following lemma isolates this device: the cutoff is opened only on the final
residual, and the remaining large prime layer gives a uniform high-frequency
Hal\'asz saving.

\begin{lemma}[Cauchy reduction and uniform high-frequency saving]\label{final_residual_pointwise}
Let $\vartheta>0$ be fixed, and let
\[
        \exp((\log X)^{3/4})\le P\le Q\le X^{o(1)}.
\]
Let
$L_*=k\ell+O(\ell^{1/3}\sqrt{\ell})$ and let
\[
        K_*=\lfloor L_*\rfloor,\qquad r=1+\ell^{-1}.
\]
Let
\[
        R_M(s)
        =
        \sum_{M<m\le2M}
        \frac{d_k(m)\1_{\Omega(m)\le L_*}}{m^s}
        \prod_{\nu=1}^{r_0}\1_{\mathcal C_\nu}(m)
        \frac1{1+\omega_{[P,Q]}(m)},
\]
where
\[
        X^\vartheta\le M\ll X .
\]
Here $r_0=O_k(1)$, and for each $\nu$ there is an interval
$I_\nu=[P_\nu,Q_\nu]\subset[2,\exp((\log X)^{3/4})]$ such that
\[
        \1_{\mathcal C_\nu}(m)
        =
        \begin{cases}
        \1_{(m,\Pi_\nu)=1},
        &\text{if }p\mid m\Rightarrow p\notin I_\nu,\\
        1-\1_{(m,\Pi_\nu)=1},
        &\text{if }\exists\,p\mid m\text{ with }p\in I_\nu,
        \end{cases}
        \qquad
        \Pi_\nu=\prod_{p\in I_\nu}p .
\]
Then there is a constant $c_0=c_0(k,\beta_0,\vartheta)>0$ such that,
uniformly for $X^{\beta_0}\le |t|\le X$,
\begin{equation}\label{pretentious_away_bound}
        R_M(1+it)
        \ll_{k,\beta_0,\vartheta}
        (\log X)^{k-1-c_0}.
\end{equation}
\end{lemma}

\begin{proof}
Since $\Omega(m)$ is integer-valued, $\Omega(m)\le L_*$ is the same condition
as $\Omega(m)\le K_*$.  Applying Lemma \ref{cauchy_cutoff_formula} with
$K=K_*$, $\rho=r$, and $a=\Omega(m)$ gives
\[
        \1_{\Omega(m)\le L_*}
        =
        \frac1{2\pi i}\int_{|z|=r}
        \frac{z^{K_*}z^{-\Omega(m)}}{z-1}\,dz .
\]
Substituting this identity into $R_M(s)$ gives
\begin{equation}\label{cauchy_Omega_cutoff}
        R_M(s)
        =
        \frac1{2\pi i}\int_{|z|=r}
        \frac{z^{K_*}}{z-1}R_M(s;z)\,dz,
\end{equation}
where at this stage $R_M(s;z)$ denotes the corresponding residual sum with the
factor $(1+\omega_{[P,Q]}(m))^{-1}$ still present.  The same lemma gives
\begin{equation}\label{cauchy_l_loss}
        \int_{|z|=r}\left|\frac{z^{K_*}}{z-1}\right|\,|dz|
        \ll_k \log\ell .
\end{equation}

For the remaining weights we use
\[
        \frac1{1+\omega_{[P,Q]}(m)}
        =\int_0^1 u^{\omega_{[P,Q]}(m)}\,du .
\]
Let
\[
        \Omega_{[P,Q]}(m)=\sum_{\substack{p^a\Vert m\\ P\le p\le Q}}a .
\]
We shall replace $u^{\omega_{[P,Q]}(m)}$ by
$u^{\Omega_{[P,Q]}(m)}$.  The two weights differ only if $p^2\mid m$ for some
$p\in[P,Q]$.  Hence, uniformly for $\Re s=1$ and $|z|=r$,
\[
\begin{aligned}
&\int_0^1
\sum_{\substack{M<m\le2M\\
u^{\omega_{[P,Q]}(m)}\ne u^{\Omega_{[P,Q]}(m)}}}
\frac{d_k(m)|z|^{-\Omega(m)}}{m}\,du  \\
&\qquad\ll_k
        \sum_{P\le p\le Q}
        \sum_{\substack{M<m\le2M\\p^2\mid m}}\frac{d_k(m)}m
        \ll_k
        (\log X)^k\sum_{p\ge P}\frac1{p^2}
        \ll_k \frac{(\log X)^k}{P}=o(1).
\end{aligned}
\]
This error is much smaller than the required bound, so it suffices to estimate
the modified integrals with $u^{\Omega_{[P,Q]}(m)}$ in place of
$u^{\omega_{[P,Q]}(m)}$.  We keep the notation $R_M(s;z)$ for this modified
quantity; the discarded error contributes $o(1)$ to \eqref{cauchy_Omega_cutoff}.
Let $\mathcal A$ be the set of indices for which the first case in the
definition of $\1_{\mathcal C_\nu}$ occurs, and let $\mathcal B$ be the set of
indices for which the second case occurs.
Also
\begin{equation}\label{local_condition_expansion}
        \prod_{\nu=1}^{r_0}\1_{\mathcal C_\nu}(m)
        =
        \prod_{\nu\in\mathcal A}\1_{(m,\Pi_\nu)=1}
        \prod_{\nu\in\mathcal B}\left(1-\1_{(m,\Pi_\nu)=1}\right)
        =
        \sum_{T\subseteq\mathcal B}(-1)^{|T|}
        \prod_{\nu\in\mathcal A\cup T}\1_{(m,\Pi_\nu)=1}.
\end{equation}
The empty product is interpreted as $1$, so the case with no remaining local
restriction is included.  Combining this expansion with the $u$-integral gives
\begin{equation}\label{R_expanded_integrands}
        R_M(s;z)
        =
        \sum_{T\subseteq\mathcal B}(-1)^{|T|}
        \int_0^1
        \sum_{M<m\le2M}\frac{g_{T,u,z}(m)}{m^s}\,du,
\end{equation}
where
\begin{equation}\label{g_explicit_definition}
        g_{T,u,z}(m)
        =
        d_k(m)z^{-\Omega(m)}u^{\Omega_{[P,Q]}(m)}
        \prod_{\nu\in\mathcal A\cup T}\1_{(m,\Pi_\nu)=1}.
\end{equation}
For fixed $T$, $u$, and $z$, let
\[
        \eta_{T,u}(p)
        =
        u^{\1_{P\le p\le Q}}
        \prod_{\nu\in\mathcal A\cup T}\1_{p\notin I_\nu},
        \qquad |\eta_{T,u}(p)|\le1.
\]
Then $g_{T,u,z}$ is multiplicative and
\[
        |g_{T,u,z}(n)|\le d_k(n),\qquad
        g_{T,u,z}(p)=kz^{-1}\eta_{T,u}(p),\qquad
        g_{T,u,z}(p^a)\ll_k a^{O_k(1)}\quad(a\ge2).
\]
Let
\[
        Y=\exp((\log X)^{3/4}).
\]
By \cite[Lemma~2]{MR15LiouvilleNote}, applied with $X$ replaced by a constant
multiple of $X$, for every fixed $C>0$ and uniformly for
$Y\le A<B\ll X$ and $X^{\beta_0}\le |t|\le X$,
\begin{equation}\label{oscillatory_prime_interval}
        \sum_{A<p\le B}\frac{p^{-it}}p
        \ll_{C,\beta_0}(\log X)^{-C}.
\end{equation}
Since all intervals $I_\nu$ lie below $Y$, the function $\eta_{T,u}$ is constant
on $O_k(1)$ prime intervals inside $(Y,M]$.  Thus \eqref{oscillatory_prime_interval}
and Mertens' estimate give
\begin{equation}\label{theta_prime_sum}
        \left|\sum_{p\le M}\eta_{T,u}(p)\frac{p^{-it}}p\right|
        \le \sum_{p\le Y}\frac1p+O_{k,C,\beta_0}((\log X)^{-C})
        \le \frac34\ell+O_k(1).
\end{equation}

Define
\[
        \mathbb D_{T,u,z}(t;M)^2
        :=
        \sum_{p\le M}
        \frac{k-\Re(g_{T,u,z}(p)p^{-it})}{p}.
\]
Using $g_{T,u,z}(p)=kz^{-1}\eta_{T,u}(p)$ and \eqref{theta_prime_sum},
\begin{align}\label{pretentious_distance_lower}
        \mathbb D_{T,u,z}(t;M)^2
        &=k\sum_{p\le M}\frac1p
          -k\Re\left(z^{-1}
          \sum_{p\le M}\eta_{T,u}(p)\frac{p^{-it}}p\right)\notag\\
        &\ge k\log\log M-\frac{3k}{4}\ell-O_{k,\beta_0}(1).
\end{align}
Since $X^\vartheta\le M\ll X$,
\[
        \log\log M=\ell+O_\vartheta(1).
\]
Consequently
\begin{equation}\label{distance_final_lower}
        \mathbb D_{T,u,z}(t;M)^2\ge \frac{k}{5}\ell
\end{equation}
for all sufficiently large $X$, uniformly in $z$, $u$, and $T\subseteq\mathcal B$.

Now the Euler factors of $g_{T,u,z}$ have the clean form
\[
        1
        \quad\text{or}\quad
        (1-\lambda_p z^{-1}p^{-s})^{-k},
        \qquad 0\le \lambda_p\le1,
\]
where $\lambda_p=\eta_{T,u}(p)$.  Hence the generalized von Mangoldt
coefficients satisfy
\[
        \Lambda_{g_{T,u,z}}(p^a)
        =k(\lambda_p z^{-1})^a\log p,
        \qquad
        |\Lambda_{g_{T,u,z}}(p^a)|\le k\log p=k\Lambda(p^a).
\]
Thus $g_{T,u,z}$ belongs uniformly to the class $\mathcal C(k)$ of
Granville--Harper--Soundararajan.  Applying
\cite[Theorem~1.1 and Corollary~1.2]{GHS19} to
$m\mapsto g_{T,u,z}(m)m^{-it}$, and using
\eqref{distance_final_lower}, gives, uniformly for $M\le y\le2M$,
\[
        \sum_{m\le y}g_{T,u,z}(m)m^{-it}
        \ll_k
        y\left\{(\log X)^{k-1}
        \exp(-c_k\ell)
        +\frac{(\log\log X)^k}{\log X}\right\},
\]
where $c_k>0$ depends only on $k$.  Partial summation on $M<m\le2M$ gives
\begin{equation}\label{halasz_input_formula}
        \sum_{M<m\le2M}\frac{g_{T,u,z}(m)}{m^{1+it}}
        \ll_k
        (\log X)^{k-1}
        \exp(-c_k\ell)
        +\frac{(\log\log X)^k}{\log X}.
\end{equation}
This gives
\[
        \sum_{M<m\le2M}\frac{g_{T,u,z}(m)}{m^{1+it}}
        \ll_{k,\beta_0,\vartheta}
        (\log X)^{k-1-c_0}
\]
with some $c_0=c_0(k,\beta_0,\vartheta)>0$.  Finally, the bounded finite
expansion, the bounded $u$-integral, and \eqref{R_expanded_integrands} give
the same bound for $R_M(1+it;z)$, uniformly for $|z|=r$.  Applying
\eqref{cauchy_Omega_cutoff} and using the Cauchy loss
\eqref{cauchy_l_loss} contributes only $(\log X)^{o(1)}$, which is absorbed by
decreasing $c_0$ slightly.  This proves \eqref{pretentious_away_bound}.
\end{proof}

\begin{proof}[Proof of Proposition \ref{MR_bound_flexible}]
Let $L$, $B$, $H$, $\mathcal S$, and $h$ be as in Proposition
\ref{MR_bound_flexible}, and set
\[
        F(s)=\sum_{\substack{X\leq n\leq4X\\ n\in\mathcal S}}\frac{d_k(n)}{n^s},
        \qquad
        T_0=\frac{X}{h}\ell.
\]
Then $Q_1T_0H/X=o(1)$.
Apply Proposition \ref{MR_modify} with
$a_n=d_k(n)\1_{\mathcal S}(n)$ and $y\in\{h,h_1\}$.  For the dyadic term in
\eqref{general_pas}, if $T>T_0$, then
\[
        \frac{X}{yT}\int_{T\leq |t|\leq2T}|F(1+it)|^2dt
        \ll \left(\frac{H}y+\frac{X}{yT}\right)(\log X)^{2k-2}
        =o((\log X)^{2k-2}).
\]
Thus it remains to focus on the range $\delta_0^{-1}\le |t|\le T_0$.
We shall prove
\begin{equation}\label{central_ms_estimate}
        \int_{\delta_0^{-1}\leq |t|\leq T_0}|F(1+it)|^2dt
        =o((\log X)^{2k-2}).
\end{equation}
On $B$ we have $d_k(n)\le k^L=H(\log X)^{k-1}$, so
\eqref{diagonal_length_property} holds with the present value of
$H$.  Let
\[
        K=\left\lfloor P_1^{1/6}(\log Q_1)^{1/3}\right\rfloor .
\]

It suffices to prove the estimate for either dyadic part of $F$.  Write the
current dyadic support as $N<n\le2N$, $N\in\{X,2X\}$, and keep the notation
$F$ for this part.  For $i=1,2$ define
\[
        \mathcal I_i=\{v\in\mathbb Z:K\log P_i\le v\le K\log Q_i\},
\]
\[
        \omega_i(m)=\#\{p\in\Primes:P_i\le p\le Q_i,\ p\mid m\},
\]
\[
        Q_{i,v}(s)=
        \sum_{\substack{P_i\le p\le Q_i\\ \exp(v/K)<p\le \exp((v+1)/K)}}
        \frac{k}{p^s}
        \qquad (v\in\mathcal I_i).
\]
For $j=1,2$ let
\[
        \mathcal A^{(j)}
        =\{m\in\mathbb N:\omega_i(m)\ge1\text{ for }i\in\{1,2\}\setminus\{j\},
        \text{ and }p\mid m\text{ for some }p\in\Primes\cap[P_3,Q_3]\}.
\]
The corresponding residual polynomial is
\[
        R_{j,v}(s)=
        \sum_{N\exp(-v/K)<m\le2N\exp(-v/K)}
        \frac{d_k(m)\1_{\Omega(m)+1\le L}\1_{\mathcal A^{(j)}}(m)}
             {m^s(1+\omega_j(m))}.
\]
The coefficients of $R_{j,v}$ are bounded by $d_k(m)\1_B(m)$.
The numbers of cells satisfy
\[
        |\mathcal I_j|\ll K\log Q_j\qquad (j=1,2).
\]
Let
\[
        \alpha_j=\frac14-\frac1{100}\left(1+\frac1{2j}\right)\qquad (j=1,2).
\]
Then $1/6<\alpha_1<\alpha_2<1/4$.
Let $\mathcal R=\{t:\delta_0^{-1}\le |t|\le T_0\}$ and let
\[
\mathcal T_1=\{t\in\mathcal R:
        |Q_{1,v}(1+it)|\le \exp(-\alpha_1 v/K)\text{ for all }v\in\mathcal I_1\},
\]
\[
\mathcal T_2=\{t\in\mathcal R\setminus\mathcal T_1:
        |Q_{2,v}(1+it)|\le \exp(-\alpha_2 v/K)\text{ for all }v\in\mathcal I_2\},
\]
and
\[
        \mathcal U=\mathcal R\setminus(\mathcal T_1\cup\mathcal T_2).
\]
Let $E_1,E_2,E_{\mathcal U}$ be the corresponding integrals of
$|F(1+it)|^2$.

\textbf{The set $\mathcal T_1$.}  Lemma \ref{mean_value_B}, applied to
$R_{1,v}$ at height $M\asymp X\exp(-v/K)$ and length $\asymp M$, gives
\begin{align}\label{E1_bound_general_new}
E_1
&\ll K\log Q_1
\sum_{v\in\mathcal I_1}
\int_{\mathcal T_1}|Q_{1,v}(1+it)R_{1,v}(1+it)|^2dt
        +o((\log X)^{2k-2})\notag\\
&\ll K\log Q_1
\sum_{v\in\mathcal I_1}\exp(-2\alpha_1v/K)
\left(1+\exp(v/K)\frac{T_0}{X}H\right)\notag\\
&\qquad\qquad\times(\log X)^{2k-2}
        +o((\log X)^{2k-2})\notag\\
&\ll P_1^{1/3-2\alpha_1+o(1)}
\left(1+Q_1\frac{T_0}{X}H\right)(\log X)^{2k-2}
        =o((\log X)^{2k-2}).
\end{align}
Here $2\alpha_1>1/3$ and $Q_1T_0H/X=o(1)$ have been used.

\textbf{The set $\mathcal T_2$.}  Since $\mathcal T_2$ is disjoint from
$\mathcal T_1$, for each $t\in\mathcal T_2$ at least one first-block cell is
large.  We split
\[
        \mathcal T_2=\bigcup_{u\in\mathcal I_1}\mathcal T_{2,u},
\]
where
\[
        \mathcal T_{2,u}:=\left\{t\in\mathcal T_2:
        u\text{ is the least cell with }
        |Q_{1,u}(1+it)|>\exp(-\alpha_1u/K)\right\}.
\]
For $u\in\mathcal I_1$ and $v\in\mathcal I_2$ let
\[
        J(u,v)=\left\lceil\frac{v/K}{u/K}\right\rceil
        \ll \frac{v/K}{\log P_1}.
\]
Applying Lemma \ref{Ramare_decop} to the second block on the partition
$\mathcal T_2=\bigcup_{u\in\mathcal I_1}\mathcal T_{2,u}$ first gives
\begin{align}
E_2
&\ll K\log Q_2
\sum_{u\in\mathcal I_1}
\sum_{v\in\mathcal I_2}
\int_{\mathcal T_{2,u}}
|Q_{2,v}(1+it)R_{2,v}(1+it)|^2dt
        +o((\log X)^{2k-2})\label{E2_Ramare_step}\\
&\ll K\log Q_2
\sum_{u\in\mathcal I_1}\sum_{v\in\mathcal I_2}
\exp\left(-\frac{2\alpha_2v}{K}
          +\frac{2\alpha_1J(u,v)u}{K}\right)\notag\\
&\qquad\qquad\times
\int_{-T_0}^{T_0}|Q_{1,u}(1+it)^{J(u,v)}R_{2,v}(1+it)|^2dt
        +o((\log X)^{2k-2}).\label{E2_amplified_step}
\end{align}

By Lemma \ref{aplf},
\begin{align}\label{E2_aplf_step}
&\int_{-T_0}^{T_0}|Q_{1,u}(1+it)^{J(u,v)}R_{2,v}(1+it)|^2dt\notag\\
&\qquad\ll
k^{2J(u,v)}(J(u,v)+1)!^2
J(u,v)\exp(u/K)(\log X)^{2k-2}.
\end{align}

We first estimate the coefficient loss in \eqref{E2_aplf_step}.  Taking
logarithms and using $J(u,v)\ll (v/K)/\log P_1$ gives
\begin{align*}
&2J(u,v)\log k+2\log((J(u,v)+1)!)
+\frac{u}{K}+\log J(u,v)\\
&\qquad\ll
\frac{v}{K}\frac{\log\log Q_2}{\log P_1}
+\log Q_1+\log\log Q_2\\
&\qquad=o(v/K).
\end{align*}
Together with
$J(u,v)u/K\le v/K+u/K=(1+o(1))v/K$, and since $v/K\ge\log P_2$,
\[
        \log(K\log Q_2|\mathcal I_1|)=o(v/K).
\]
The final summation over $v$ contributes only $\exp(o(\log P_2))$.  Substituting
\eqref{E2_aplf_step} into \eqref{E2_amplified_step}, we obtain
\begin{align}\label{E2_bound_general_new}
E_2
&\ll(\log X)^{2k-2}
\sum_{v\in\mathcal I_2}
\exp\left(-2(\alpha_2-\alpha_1+o(1))\frac{v}{K}\right)
        +o((\log X)^{2k-2})\notag\\
&\ll(\log X)^{2k-2}
\exp\{-2(\alpha_2-\alpha_1+o(1))\log P_2\}
        =o((\log X)^{2k-2}).
\end{align}

\textbf{The set $\mathcal U$.}  For $v\in\mathcal I_2$ let
\[
        Y_{2,v}=\exp(v/K),
        \qquad V_{2,v}=Y_{2,v}^{\alpha_2}.
\]
Here the first two prime blocks have large cells; the third block is handled by
a final Ramar\'e decomposition and a small/large split.
We follow the large-value treatment in \cite[Section~8]{MR16}.  Define
\[
        \mathcal I_3=\{w\in\mathbb Z:K\log P_3\le w\le K\log Q_3\},
        \qquad
        \omega_3(m)=\#\{p\in\Primes:P_3\le p\le Q_3,\ p\mid m\}.
\]
For $w\in\mathcal I_3$ define
\[
        Q_{3,w}(s)=
        \sum_{\substack{P_3\le p\le Q_3\\ \exp(w/K)<p\le \exp((w+1)/K)}}
        \frac{k}{p^s},
\]
and
\[
        R_{3,w}(s)=
        \sum_{\substack{N\exp(-w/K)<m\le2N\exp(-w/K)\\ \Omega(m)+1\le L}}
        \frac{d_k(m)\1_{\omega_1(m)\ge1}\1_{\omega_2(m)\ge1}}
             {m^s(1+\omega_3(m))}.
\]
Applying Lemma \ref{Ramare_decop} to the third block gives
\begin{equation}\label{EU_ramare_new}
E_{\mathcal U}
\ll K\left(1+\log\frac{Q_3}{P_3}\right)
\sum_{w\in\mathcal I_3}
        \int_{\mathcal U}
        |Q_{3,w}(1+it)R_{3,w}(1+it)|^2dt
        +o((\log X)^{2k-2}).
\end{equation}
The error term is
\[
        \ll \left(\frac1K+\frac1{P_3}\right)
        \left(\frac{T_0}{X}H+1\right)(\log X)^{2k-2}
        =o((\log X)^{2k-2}).
\]

Fix $w\in\mathcal I_3$.  Choose a well-spaced set
$\mathcal V_w\subset\mathcal U$ by taking one point from every other occupied
unit interval so that
\[
\int_{\mathcal U}|Q_{3,w}(1+it)R_{3,w}(1+it)|^2dt
\ll
\sum_{t\in\mathcal V_w}|Q_{3,w}(1+it)R_{3,w}(1+it)|^2 .
\]
Since $\mathcal U\cap\mathcal T_2=\emptyset$, every $t\in\mathcal V_w$ has a
second-block cell $v\in\mathcal I_2$ with
$|Q_{2,v}(1+it)|>V_{2,v}^{-1}$.  Put
\[
        \mathcal V_w(v)=\{t\in\mathcal V_w:
        |Q_{2,v}(1+it)|>V_{2,v}^{-1}\}.
\]
Lemma \ref{large_value_results}(ii), applied with
$P_0=Y_{2,v}$, $V=V_{2,v}=Y_{2,v}^{\alpha_2}$, and $T=T_0$, gives
\[
        |\mathcal V_w(v)|
        \ll T_0^{2\alpha_2}
        Y_{2,v}^{2\alpha_2}
        \exp\left(2k\frac{\log T_0}{\log Y_{2,v}}\log\log T_0\right).
\]
Since $\log Y_{2,v}\ge \log P_2$ and $\log P_2/\ell\to\infty$, the final
exponential is $X^{o(1)}$.  Summing over $v\in\mathcal I_2$ therefore gives
\begin{equation}\label{U_cardinality_second_block}
        |\mathcal V_w|\ll X^{2\alpha_2+o(1)}.
\end{equation}
Here we also use $T_0\le X$, $Y_{2,v}\le Q_2=X^{o(1)}$, and
$|\mathcal I_2|\ll K\log Q_2=X^{o(1)}$.  Thus $|\mathcal V_w|\ll X^{1/2-\eta_1}$ for some
fixed $\eta_1>0$.

Let
\[
        Z=(\log X)^{k^2/2+10}.
\]
Define
\[
        \mathcal V_w^S=\{t\in\mathcal V_w:
        |Q_{3,w}(1+it)|\le Z^{-1}\},
\]
and
\[
        \mathcal V_w^L=\{t\in\mathcal V_w:
        |Q_{3,w}(1+it)|> Z^{-1}\}.
\]
\begin{equation}\label{EU_small_final_cell}
K\left(1+\log\frac{Q_3}{P_3}\right)
\sum_{w\in\mathcal I_3}
\sum_{t\in\mathcal V_w^S}|Q_{3,w}(1+it)R_{3,w}(1+it)|^2
=o((\log X)^{2k-2}).
\end{equation}
Indeed, let $M_w=N\exp(-w/K)$.  Since $|Q_{3,w}(1+it)|\le Z^{-1}$
on $\mathcal V_w^S$, Lemma \ref{discrete_mean_value_theorem} gives
\begin{align*}
&\sum_{w\in\mathcal I_3}
\sum_{t\in\mathcal V_w^S}|Q_{3,w}(1+it)R_{3,w}(1+it)|^2\\
&\qquad\ll
Z^{-2}\sum_{w\in\mathcal I_3}
\sum_{t\in\mathcal V_w}|R_{3,w}(1+it)|^2\\
&\qquad\ll
Z^{-2}\sum_{w\in\mathcal I_3}
\left(1+\frac{|\mathcal V_w|T_0^{1/2}}{M_w}\right)
(\log X)
\frac1{M_w}\sum_{M_w<m\le2M_w}d_k(m)^2\1_B(m).
\end{align*}
Here $M_w=X^{1-o(1)}$ and \eqref{U_cardinality_second_block} gives
$|\mathcal V_w|\ll X^{1/2-\eta_1}$, so the factor in parentheses is
$1+o(1)$.  Also, by \eqref{diagonal_length_property} and Shiu's estimate
\eqref{Shiu_weighted} with $\eta=0$,
\[
        \frac1{M_w}\sum_{M_w<m\le2M_w}d_k(m)^2\1_B(m)
        \ll H(\log X)^{2k-2}.
\]
Thus
\[
\sum_{w\in\mathcal I_3}
\sum_{t\in\mathcal V_w^S}|Q_{3,w}(1+it)R_{3,w}(1+it)|^2
\ll
(\log X)^{-k^2-19+o(1)}
|\mathcal I_3|H(\log X)^{2k-2}.
\]
Multiplying by $K(1+\log(Q_3/P_3))$ and using
$|\mathcal I_3|\ll K\log Q_3$ gives a bound
\[
        \ll
        (\log X)^{-k^2-19+o(1)}
        H(K\log Q_3+3)^2(\log X)^{2k-2}
        =o((\log X)^{2k-2}),
\]
since $|\nu(X)|\le\ell^{1/3}$ gives
$H=(\log X)^{k\log k-k+1+o(1)}$, while
$K=(\log X)^{o(1)}$ and $\log Q_3\le\log X$.

For the large part, Lemma \ref{large_value_results}(ii) gives
\[
|\mathcal V_w^L|
\le
\exp\{O_k(\ell^3)\}.
\]
Since $w/K\ge\log P_3=(\log X)/\ell^2$, Lemma
\ref{short_prime_cell_square}, applied with $\varepsilon=1/100$, yields
\[
|\mathcal V_w^L|
\exp\left(-\frac{w/K}{(\log T_0)^{2/3+1/100}}\right)(\log T_0)^2=o(1).
\]
Consequently
\begin{equation}\label{EU_large_prime_square_sum}
\sum_{t\in\mathcal V_w^L}|Q_{3,w}(1+it)|^2
\ll \frac1{K(\log P_3)^2}.
\end{equation}
Here $R_{3,w}$ is supported on $M\asymp X\exp(-w/K)$, and
$\exp(w/K)\le Q_3=X^{o(1)}$; hence $X^{1/2}\le M\ll X$ for all large $X$.
Also $P_3\ge\exp((\log X)^{3/4})$, while the local restrictions inherited from
the first two blocks lie below $\exp((\log X)^{3/4})$.
Lemma \ref{final_residual_pointwise} gives, uniformly on this high-frequency
range,
\[
        |R_{3,w}(1+it)|^2
        \ll(\log X)^{2k-2-2c_0+o(1)}
\]
after decreasing $c_0$ slightly if necessary.  Using
\eqref{EU_large_prime_square_sum} and summing over $w$, we obtain
\begin{align}\label{EU_large_final_cell}
&K\left(1+\log\frac{Q_3}{P_3}\right)
\sum_{w\in\mathcal I_3}
\sum_{t\in\mathcal V_w^L}
|Q_{3,w}(1+it)R_{3,w}(1+it)|^2 \notag\\
&\qquad\ll
(\log X)^{2k-2-2c_0+o(1)}
\frac{\left(1+\log(Q_3/P_3)\right)|\mathcal I_3|}{(\log P_3)^2}.
\end{align}
Since $|\mathcal I_3|\ll K\log Q_3$ and $K\le(\log X)^{o(1)}$,
\[
        \frac{\left(1+\log(Q_3/P_3)\right)|\mathcal I_3|}{(\log P_3)^2}
        \ll
        K\left(\frac{\log Q_3}{\log P_3}\right)^2
        =(\log X)^{o(1)} ,
\]
and therefore the right-hand side of \eqref{EU_large_final_cell} is
\[
        \ll
        (\log X)^{2k-2-2c_0+o(1)}
        =o((\log X)^{2k-2})
\]
for the constant $c_0>0$ in Lemma \ref{final_residual_pointwise}.  Combining
\eqref{EU_ramare_new}, \eqref{EU_small_final_cell}, and
\eqref{EU_large_final_cell} gives
\begin{equation}\label{EU_sum_final}
        E_{\mathcal U}=o((\log X)^{2k-2}).
\end{equation}
Adding \eqref{E1_bound_general_new}, \eqref{E2_bound_general_new},
and \eqref{EU_sum_final} proves \eqref{central_ms_estimate}.

For the dyadic maximum in \eqref{general_pas}, the part with $2T\le T_0$ is
covered by \eqref{central_ms_estimate}, since $T\ge X/y$ implies
$X/(yT)\le1$.  If $T<T_0<2T$, split the integral at $T_0$: the lower part is
covered by \eqref{central_ms_estimate}, while the upper part is bounded by the
same mean-value estimate used above, since $T>T_0/2$.  The case $T\ge T_0$ is
also covered by that estimate.
Therefore Proposition \ref{MR_modify} gives, for $y=h$ and $y=h_1$, an $L^2$
error $o((\log X)^{2k-2})$ after replacing the short sum by its central Perron
integral.  Lemma
\ref{central_int} compares the two central Perron integrals, so
\[
\frac1X\int_X^{2X}\left|
        \frac1h\sum_{\substack{x<n\leq x+h\\ n\in\mathcal S}}d_k(n)
        -
        \frac1{h_1}\sum_{\substack{x<n\leq x+h_1\\ n\in\mathcal S}}d_k(n)
\right|^2dx=o((\log X)^{2k-2}).
\]
Cauchy's inequality gives the asserted $L^1$ estimate.  Since $h$ and $h_1$
are integers, the two normalised sums are constant on each interval $[m,m+1)$;
Markov's inequality gives the integer almost-all assertion.
\end{proof}

\bibliography{divisor_critical_window_referee_revised}

@book{IK04,
  author    = {Iwaniec, Henryk and Kowalski, Emmanuel},
  title     = {Analytic Number Theory},
  series    = {American Mathematical Society Colloquium Publications},
  volume    = {53},
  publisher = {American Mathematical Society},
  address   = {Providence, RI},
  year      = {2004}
}

@article{Henriot12,
  author  = {Henriot, Kevin},
  title   = {Nair--Tenenbaum bounds uniform with respect to the discriminant},
  journal = {Mathematical Proceedings of the Cambridge Philosophical Society},
  volume  = {152},
  number  = {3},
  pages   = {405--424},
  year    = {2012},
  doi     = {10.1017/S0305004111000752}
}

@article{Henriot14erratum,
  author  = {Henriot, Kevin},
  title   = {Nair--Tenenbaum uniform with respect to the discriminant---Erratum},
  journal = {Mathematical Proceedings of the Cambridge Philosophical Society},
  volume  = {157},
  number  = {2},
  pages   = {375--377},
  year    = {2014},
  doi     = {10.1017/S0305004114000341}
}

@article{KMS22,
  author  = {Khan, Rizwanur and Milinovich, Micah B. and Subedi, Unique},
  title   = {A weighted version of the {Erd\H{o}s}--{Kac} theorem},
  journal = {Journal of Number Theory},
  volume  = {239},
  pages   = {1--20},
  year    = {2022}
}

@article{mangerel2021divisor,
  author  = {Mangerel, Alexander P.},
  title   = {Divisor-bounded multiplicative functions in short intervals},
  journal = {Research in the Mathematical Sciences},
  volume  = {10},
  number  = {12},
  year    = {2023},
  doi     = {10.1007/s40687-023-00376-0}
}

@article{MR16,
  author  = {Matom\"aki, Kaisa and Radziwi{\l\l}, Maksym},
  title   = {Multiplicative functions in short intervals},
  journal = {Annals of Mathematics. Second Series},
  volume  = {183},
  number  = {3},
  pages   = {1015--1056},
  year    = {2016}
}

@misc{matomaki2020multiplicative,
  author        = {Matom\"aki, Kaisa and Radziwi{\l\l}, Maksym},
  title         = {Multiplicative functions in short intervals {II}},
  eprint        = {2007.04290},
  archivePrefix = {arXiv},
  primaryClass  = {math.NT},
  year          = {2020}
}

@article{MR3961326,
  author  = {Matom\"aki, Kaisa and Radziwi{\l\l}, Maksym and Tao, Terence},
  title   = {Correlations of the von {Mangoldt} and higher divisor functions {II}: divisor correlations in short ranges},
  journal = {Mathematische Annalen},
  volume  = {374},
  number  = {1-2},
  pages   = {793--840},
  year    = {2019}
}

@article{MSTT23,
  author  = {Matom\"aki, Kaisa and Shao, Xuancheng and Tao, Terence and Ter\"av\"ainen, Joni},
  title   = {Higher uniformity of arithmetic functions in short intervals. {I}. {All} intervals},
  journal = {Forum of Mathematics, Pi},
  volume  = {11},
  pages   = {Paper No. e29, 97},
  year    = {2023},
  doi     = {10.1017/fmp.2023.28}
}

@article{MRSTT26,
  author  = {Matom\"aki, Kaisa and Radziwi{\l\l}, Maksym and Shao, Xuancheng and Tao, Terence and Ter\"av\"ainen, Joni},
  title   = {Higher uniformity of arithmetic functions in short intervals. {II}. {Almost} all intervals},
  journal = {Inventiones Mathematicae},
  volume  = {244},
  pages   = {967--1091},
  year    = {2026},
  doi     = {10.1007/s00222-026-01408-6}
}

@article{Sh80,
  author  = {Shiu, P.},
  title   = {A {Brun}--{Titchmarsh} theorem for multiplicative functions},
  journal = {Journal f\"ur die Reine und Angewandte Mathematik},
  volume  = {313},
  pages   = {161--170},
  year    = {1980}
}

@article{LiuWu21,
  author  = {Liu, Kui and Wu, Jie},
  title   = {Weighted {Erd\H{o}s}--{Kac} theorem in short intervals},
  journal = {The Ramanujan Journal},
  volume  = {55},
  number  = {1},
  pages   = {1--12},
  year    = {2021},
  doi     = {10.1007/s11139-020-00343-1}
}

@book{Te15,
  author    = {Tenenbaum, G{\'e}rald},
  title     = {Introduction to Analytic and Probabilistic Number Theory},
  series    = {Graduate Studies in Mathematics},
  volume    = {163},
  publisher = {American Mathematical Society},
  address   = {Providence, RI},
  year      = {2015}
}
\bibliographystyle{plain}
\end{document}